\documentclass[preprint,12pt, a4paper]{elsarticle}

\usepackage{amssymb}
\usepackage{amsmath}
\usepackage{hyperref}
\usepackage{xurl}
\usepackage{multirow}
\usepackage{amsthm}
\usepackage{graphicx}
\usepackage{booktabs}
\usepackage{multirow}

\newtheorem{remark}{Remark}
\journal{SoftwareX}

\begin{document}
\renewcommand{\labelenumii}{\arabic{enumi}.\arabic{enumii}}

\begin{frontmatter}



\title{$\alpha\beta$-Tessarine Toolbox: A high-performance MATLAB framework for hypercomplex tensor algebra}


\author[label1]{José Domingo Jiménez-López}
\author[label1]{Jesús Navarro-Moreno}
\author[label1]{Juan Carlos Ruiz-Molina}
\address[label1]{Department of Statistics and Operations
Research, University of Ja\'en, Campus Las Lagunillas s/n, Ja\'en, 23071, Spain
, e-mails: \{jdomingo, jnavarro, jcruiz\}@ujaen.es}

\begin{abstract}

The \textit{$\alpha\beta$-Tessarine Toolbox} is introduced as a high-performance MATLAB framework for generalized hypercomplex tensor algebra. By implementing the $\alpha\beta$-tessarine mathematical structure, a strictly more general environment than existing quaternion or limited tessarine packages is provided. This open-source software bridges abstract theoretical models and experimental applications, allowing researchers to efficiently solve multidimensional problems. Built-in capabilities include exact hypercomplex matrix factorizations, such as Singular Value Decomposition (SVD) and QR, optimized for large-scale image in-painting, denoising, and digital watermarking. The scalable architecture facilitates new computational research in machine learning and multidimensional signal processing.

\end{abstract}

\begin{keyword}
$\alpha\beta$-tessarines \sep Hypercomplex algebra \sep Tensor processing \sep MATLAB toolbox \sep Multidimensional signals \sep High-performance computing

\end{keyword}

\end{frontmatter}


\section*{Required Metadata}
\label{}

\section*{Current code version}
\label{}

\begin{table}[!h]
\begin{tabular}{|l|p{6.5cm}|p{6.5cm}|}
\hline
\textbf{Nr.} & \textbf{Code metadata description} & \textbf{Metadata} \\
\hline
C1 & Current code version & v1.0.0 \\
\hline
C2 & Permanent GitHub link to code/repository used for this code version & \url{https://github.com/jdomingo2003/abtessarine-toolbox}\\
\hline
C3 & Legal Code License   & GNU General Public License v3.0 (GPL-3.0) \\
\hline
C4 & Code versioning system used & Git  \\
\hline
C5 & Software code languages, tools, and services used & MATLAB \\
\hline
C6 & Compilation requirements, operating environments \& dependencies & MATLAB (R2017A or newer). No compilation required. Cross-platform (Windows, macOS, Linux).  \\
\hline
C7 & If available Link to developer documentation/manual & \url{https://github.com/jdomingo2003/abtessarine-toolbox/blob/main/README.md}\\
\hline
C8 & Support email for questions & jdomingo@ujaen.es; jnavarro@ujaen.es; jcruiz@ujaen.es\\
\hline
\end{tabular}
\caption{Code metadata}
\label{code_metadata} 
\end{table}

\section{Motivation and significance}
Multidimensional phenomena in modern physics, engineering, and data science are increasingly modeled using hypercomplex number systems \cite{Hahn2016, Valle2021, Vieira2022, ZengSong, Borio, ZhangGao}. Historically, frameworks derived from the Cayley-Dickson construction (spanning from standard quaternions to higher-dimensional algebras like sedenions and trigintaduonions) have dominated the literature. However, exploiting these algebras computationally presents a major hurdle: as dimensions increase, fundamental algebraic properties such as commutativity and associativity are irreversibly lost. This structural deficiency severely complicates standard matrix operations and spectral theory.

To circumvent the limitations of non-commutative spaces, researchers have pivoted towards associative and commutative 4-dimensional (4D) hypercomplex systems, most notably tessarines \cite{Pei, Ortolani2017, Cariow2024}. The mathematical flexibility of tessarines has recently inspired the formulation of parameterized algebraic families, including $\beta$-quaternions \cite{Navarro_beta_quaternion} and generalized Segre’s quaternions (GSQ) \cite{Navarro_Segre}. These generalized spaces effectively govern specific mathematical branches, such as elliptic  quaternions \cite{Kosal1, Kosal3}, and hyperbolic quaternions, depending on their topological parameters.

Despite the theoretical breakthroughs in parameterized commutative algebras, their transition into applied sciences is critically hindered by a software gap. The lack of standardized computational tools forces researchers to develop ad-hoc, inefficient code. The history of hypercomplex research demonstrates that algorithmic accessibility drives scientific discovery; for instance, the widespread adoption of quaternion algebra was fundamentally catalyzed by Sangwine and Bihan’s MATLAB toolbox in 2005 \cite{Sanqwine_LeBihan}, and similar early-stage efforts have recently appeared for specific elliptic matrices \cite{Kosal2}. Recently, a highly significant advancement in this mathematical landscape was presented in \cite{JimenezLopez2025}, where the novel family of commutative $\alpha\beta$-tessarine algebras was formally defined. This algebraic system strictly generalises Segre’s and elliptic quaternions, demonstrating profound computational advantages over traditional non-commutative quaternions by maintaining compatibility with standard linear algebra techniques. The authors in [\cite{JimenezLopez2025} established a rigorous theoretical foundation for $\alpha\beta$-tessarine matrices, formulating exact mechanisms for inversion, LU factorization with partial pivoting, determinant calculus, and a comprehensive spectral theory encompassing Singular Value Decomposition (SVD), rank-$k$ approximations, and pseudoinverses. However, despite its proven potential to outperform traditional frameworks in real-world scenarios such as image reconstruction and face recognition, this powerful algebra entirely lacks a dedicated, open-source computational infrastructure accessible to the broader scientific community. Without a native computational framework, implementing advanced statistical techniques—such as exploiting properness to reduce the effective dimensionality of stochastic processes \cite{Bihan2017, Nitta2019, Grassucci_proper}—remains computationally prohibitive.

The \texttt{abtessarine\_Toolbox} is specifically engineered to bridge this gap by delivering the first robust computational environment for the $\alpha\beta$-tessarine algebra and its 8-dimensional (8D) extension, the generalized $\alpha\beta$-tessarine algebra. By fully supporting an algebra that encompasses both GSQ and elliptic quaternions, this software empowers the scientific community to deploy highly complex spectral analyses and translates the theoretical breakthroughs of \cite{JimenezLopez2025} into an immediate, deployable technological reality. 

This toolbox will directly accelerate research in augmented statistics, optimal control, and non-circular signal processing. It provides mathematically exact and computationally optimized implementations of essential operations, ranging from basic algebraic inverses and square roots to advanced matrix factorizations (LU, SVD, Cholesky) and rank-k approximations. Consequently, researchers can now validate theoretical frameworks and develop novel adaptive filters \cite{Navarro_beta_quaternion} or deep learning algorithms, among other techniques, directly in the $\alpha\beta$-tessarine domain without encountering mathematical contradictions. Much like Sangwine and Le Bihan’s toolbox catalyzed decades of quaternion research, this software aims to serve as the foundational catalyst for a prolific new era of experimental research across diverse scientific disciplines.

The toolbox architecture leverages MATLAB's Object-Oriented Programming (OOP) to provide a transparent, zero-overhead user experience. The workflow is highly intuitive:

\begin{enumerate}
    \item \textit{Topological Setup:} The user globally defines the algebraic space by initializing the scalar parameters ($\alpha$ and $\beta$) via the \texttt{setabtessarine} function.
    \item \textit{Matrix Generation:} Hypercomplex arrays are instantiated using custom constructors (e.g., \texttt{abtzeros}, \texttt{abteye}, \texttt{abtrandn}), which handle the internal memory allocation of the real and complex branches.
    \item \textit{Overloaded Execution:} Operations are executed using native MATLAB syntax. Standard mathematical symbols (e.g., \texttt{+}, \texttt{*}, \texttt{/}) and functions (\texttt{svd(X)}, \texttt{eig(X)}, \texttt{inv(X)}) are deeply overloaded. This allows the user to process hypercomplex arrays exactly as if they were built-in numeric types, while the software safely manages the underlying isomorphic transformations.
\end{enumerate}

\section{Software description}

The \texttt{abtessarine\_Toolbox} is a MATLAB-based computational framework designed to operationalize the theoretical foundations of commutative hypercomplex algebras, specifically the $\alpha\beta$-tessarines and their 8D generalizations. The software provides a robust environment to execute advanced linear algebra and spectral operations. In addition to directly translating the mathematical theorems established by Jim{\'{e}}nez-L{\'{o}}pez et al. \cite{JimenezLopez2025} into computationally efficient algorithms, the toolbox introduces the mathematical formulation and implementation of the hypercomplex Cholesky factorization, whose exact theoretical development is provided in Appendix A.

\subsection{Software architecture}

To guarantee a seamless user experience and computational efficiency, the software is built entirely upon MATLAB's OOP paradigm. The architecture avoids computationally expensive nested \texttt{for}-loops by representing hypercomplex numbers through isomorphic matrices based on real and complex branches.

The core architecture is divided into two primary classes:
\begin{itemize}
    \item \texttt{@abtessarine}: Handles 4D $\alpha\beta$-tessarine objects.
    \item \texttt{@gabtessarine}: Manages the 8D generalized $\alpha\beta$-tessarine objects.
\end{itemize}

Internally, when an object is instantiated, the toolbox leverages the globally defined topological parameters ($\alpha$ and $\beta$) to map the hypercomplex arrays into their equivalent real or complex block-matrix representations. This architectural decision is crucial, as it allows the toolbox to exploit MATLAB's highly optimized, multi-threaded Basic Linear Algebra Subprograms (BLAS) and Linear Algebra PACKage (LAPACK) backends for heavy computations, before mapping the results back to the hypercomplex domain.

\begin{figure}[h]
    \centering
    \includegraphics[width=0.8\textwidth]{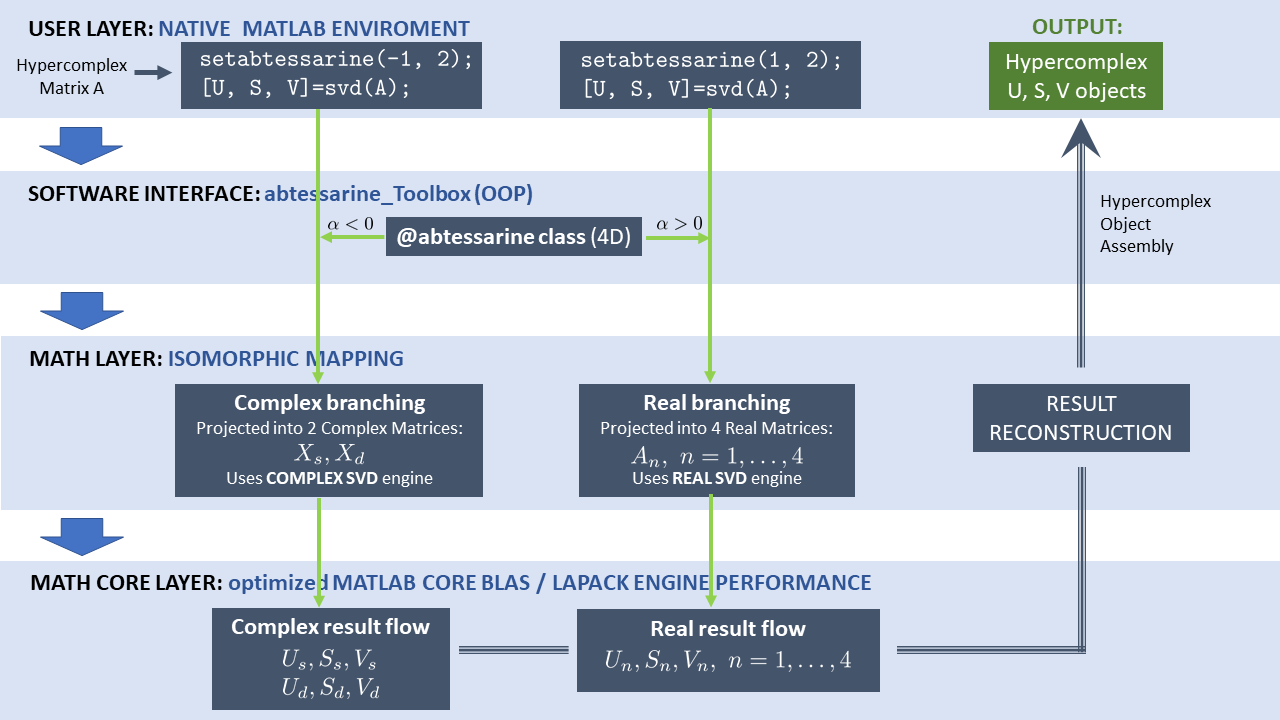}
    \caption{Architectural workflow of the \texttt{abtessarine\_Toolbox}, illustrating the data flow from native MATLAB user input, through the isomorphic mapping layer, to the core BLAS/LAPACK execution engine.}
    \label{fig:architecture}
\end{figure}

\subsection{Software functionalities}

The primary goal of the toolbox is to provide a transparent interface where users can manipulate hypercomplex matrices exactly as they would with standard numeric arrays. As detailed in the theoretical framework \cite{JimenezLopez2025}, the functionalities implemented are strictly mathematically consistent with the properties of commutative hypercomplex algebras. 

The major functionalities of the software include:
\begin{itemize}
    \item \textit{Environment Initialization and Memory Allocation:} Functions such as \texttt{setabtessarine} establish the global topological space. Allocators like \texttt{abtzeros}, \texttt{abtones}, \texttt{abteye}, and \texttt{abtrandn} are provided for rapid matrix generation.
    \item \textit{Basic Arithmetic and Logical Operations:} Deep overloading of standard operators (\texttt{+}, \texttt{-}, \texttt{*}, \texttt{.*}, \texttt{/}) ensuring the preservation of commutativity and associativity within the hypercomplex domain.
    \item \textit{Matrix Factorizations:} Native execution of complex matrix decompositions, including LU factorization with partial pivoting (\texttt{lu}), QR decomposition (\texttt{qr}), and Cholesky factorization (\texttt{chol}) for positive-definite hypercomplex matrices.
    \item \textit{Spectral Theory and Inversion:} Implementation of eigenvalues and eigenvectors computation (\texttt{eig}), SVD (\texttt{svd}), mathematical inverses (\texttt{inv}), and Moore-Penrose pseudoinverses (\texttt{pinv}) using the exact methods formulated in \cite{JimenezLopez2025}.
    \item \textit{Optimization and Least Squares:} Direct integration of the backslash operator (\texttt{\textbackslash}) to solve minimum norm least-squares problems for overdetermined systems.
\end{itemize}

A comprehensive set of sample code snippets, detailing the initialization of the algebra, memory-efficient generation, and the execution of 4D and 8D advanced spectral decompositions, is provided in the well-documented \href{https://github.com/jdomingo2003/abtessarine-toolbox/blob/main/README.md}{\texttt{README.md}} file within the software's official GitHub repository.

To guarantee the mathematical integrity of all overloaded operators and ensure the long-term maintainability of the software, a comprehensive validation framework has been integrated into the repository's \href{https://github.com/jdomingo2003/abtessarine-toolbox/tree/main/tests}{\texttt{tests/}} directory. This framework includes both high-level algebraic stress tests and unit-level verification scripts tailored for every 4D and 8D hypercomplex method. By executing these test suites, researchers can independently verify code correctness and ensure strict numerical precision across different hardware architectures, establishing a reliable foundation for future extensions of the toolbox.

\section{Illustrative examples}

To rigorously evaluate the computational efficiency, mathematical robustness, and physical scalability of the proposed \texttt{abtessarine\_Toolbox}, three diverse, high-performance signal processing applications were designed. In all scenarios, the $\alpha\beta$-tessarine techniques are directly benchmarked against their quaternion counterparts under strictly identical conditions. Crucially, to ensure a highly competitive and realistic comparison, the quaternion baseline does not rely on the legacy Sangwine and Le Bihan's toolbox \cite{Sanqwine_LeBihan}. Instead, it employs heavily optimized, custom-vectorized quaternion functions that significantly outperform traditional implementations. This guarantees that the observed performance gains are intrinsically derived from the algebraic and structural advantages of the commutative framework, rather than mere disparities in software optimization. The selected applications are:
\begin{enumerate}
    \item \textit{Color Image In-painting via Hypercomplex Singular Value Thresholding (SVT) (Section~\ref{sec:inpainting}):} An iterative, mathematically ill-posed inverse problem utilizing subspace-based truncated SVD solvers (\texttt{svds}) to recover missing spatial data under varying damage saturation levels. 
    [\href{https://github.com/jdomingo2003/abtessarine-toolbox/blob/main/examples/example11_inpainting_San_Francisco.m}{\texttt{Code Exp I}} | \href{https://github.com/jdomingo2003/abtessarine-toolbox/blob/main/examples/example12_inpainting_Autumn_Forest.m}{\texttt{Code Exp II}}]
    
    \item \textit{Global Image Denoising via Dense SVD (Section~\ref{sec:denoising}):} A highly demanding, compute-bound stress test employing full, dense SVD decompositions to filter out additive white Gaussian noise (AWGN) across progressively scaling square spatial resolutions. 
    [\href{https://github.com/jdomingo2003/abtessarine-toolbox/blob/main/examples/example2_denoising.m}{\texttt{Code}}]
    
    \item \textit{Digital Image Watermarking via Hypercomplex QR Decomposition (Section~\ref{sec:watermarking}):} An energy-preserving security application developed to evaluate algebraic orthogonality, numerical precision, and $\mathcal{O}(N^3)$ computational complexity under strict isometric constraints. 
    [\href{https://github.com/jdomingo2003/abtessarine-toolbox/blob/main/examples/example3_watermarked.m}{\texttt{Code}}]
\end{enumerate}

To ensure baseline reproducibility and provide a clear computational context, all numerical experiments were performed on a standard laptop equipped with a 13th Gen Intel\textsuperscript{\textregistered} Core\textsuperscript{\texttrademark} i9-13980HX 2.20 GHz processor and 32 GB of RAM.

For a scientifically rigorous evaluation, a heterogeneous dataset comprising high-resolution and ultra-high-resolution images was selected. Rather than relying on standard low-resolution academic examples, massive spatial dimensions were deliberately targeted by the benchmarking pipeline. This approach was chosen to thoroughly expose the memory boundaries and computational bottlenecks of hypercomplex representations under realistic, extreme-scale conditions. 

Detailed specifications for the selected dataset, including exact pixel dimensions, licensing details, and direct source links to guarantee full reproducibility, are summarized in Table~\ref{tab:images}. Furthermore, to facilitate immediate access and testing, the source files for these images have been included directly within the repository's \href{https://github.com/jdomingo2003/abtessarine-toolbox/tree/main/examples}{\texttt{examples/}} folder. The only exception is the ultra-high-resolution \textit{Carina Nebula} image, which, due to GitHub's strict file size limitations, could not be hosted natively in the repository; however, it remains fully accessible for direct download via its corresponding source link in Table \ref{tab:images}.
\begin{remark}
The inclusion of the ultra-high-resolution \textit{Carina Nebula} dataset (comprising approximately 123 megapixels) highlights a critical architectural advantage of the proposed commutative $\alpha\beta$-tessarine framework over existing quaternion baselines, particularly regarding memory footprint and computational bottlenecks. In a traditional quaternion-based signal processing pipeline, the complex isomorphic mapping required to perform spectral decompositions would force the construction of a complex matrix of size $29150 \times 16882$. Stored in a double-precision complex format, this isomorphic matrix would occupy approximately 7.33 GB of physical RAM solely for input storage, frequently exceeding the heap limits of standard commodity workstations.

In contrast, by exploiting the idempotent representation of commutative $\alpha\beta$-tessarines, the problem is decoupled by the proposed toolbox into independent isomorphic representations of size $14575 \times 8441$. This fundamental decoupling strictly halves the peak memory requirements and avoids the severe latency penalties associated with caching and manipulating massively scaled dense matrices. As will be empirically demonstrated in subsequent SVD benchmarks, this structural efficiency not only prevents memory overflows but also enables dramatic reductions in execution time.
\end{remark}

\begin{table*}[htbp]
\centering
\caption{Benchmark image specifications, native dimensions, licensing models, and direct HTTPS source links.}
\label{tab:images}
\renewcommand{\arraystretch}{1.2} 
\resizebox{\textwidth}{!}{%
\begin{tabular}{@{}lllc@{}} 
\toprule
\textbf{Image Name} & \textbf{Dimensions} & \textbf{Credit \& License} & \textbf{Source Link} \\ \midrule
Autumn Forest & $3640 \times 5464$ & Unsplash License (ML6kHR--Uys) & \href{https://unsplash.com/es/fotos/un-bosque-lleno-de-muchos-arboles-de-diferentes-colores-ML6kHR--Uys}{\texttt{[Link]}} \\ \addlinespace
Modern Architecture & $3911 \times 5867$ & Unsplash License (Opwvoz9zwYk) & \href{https://unsplash.com/es/fotos/un-edificio-muy-alto-con-muchas-ventanas-Opwvoz9zwYk}{\texttt{[Link]}} \\ \addlinespace
San Francisco & $7952 \times 5304$ & Unsplash License (o8Utw2ETExA) & \href{https://unsplash.com/es/fotos/puente-golden-state-san-francisco-o8Utw2ETExA}{\texttt{[Link]}} \\ \addlinespace
Carina Nebula & $14575 \times 8441$ & NASA/ESA/CSA/STScI (Public Domain) & \href{https://science.nasa.gov/asset/webb/cosmic-cliffs-in-the-carina-nebula-nircam-image/}{\texttt{[Link]}} \\ \bottomrule
\end{tabular}%
}
\end{table*}

\subsection{Illustrative Example: Color Image In-painting via Hypercomplex SVT}
\label{sec:inpainting}

In this subsection, the practical applicability and computational superiority of the \texttt{abtessarine\_Toolbox} are demonstrated by solving a large-scale inverse problem: color image in-painting. Through the implementation of an iterative hypercomplex SVT algorithm, missing visual data is reconstructed across massive high-resolution datasets. The following experiments rigorously benchmark the proposed commutative framework against highly optimized quaternion baselines, evaluating restoration fidelity by the Peak Signal-to-Noise Ratio (PSNR), execution speedups, and the critical ability of the toolbox to bypass severe memory and algorithmic bottlenecks under extreme dimensional constraints.

\subsubsection{Mathematical Foundations of Hypercomplex In-painting}
Digital color image in-painting is a mathematically ill-posed inverse problem aimed at restoring missing, corrupted, or occluded pixels within a visual scene. In this framework, a color image is modeled as a hypercomplex matrix $Z \in \mathbb{A}^{M \times N}$ (where $\mathbb{A}$ denotes either the $\alpha\beta$-tessarine or quaternion algebra), mapping the RGB color channels to the imaginary components of the hypercomplex structure. The underlying assumption enabling this reconstruction is that natural scenes exhibit a high degree of non-local self-similarity and redundant structural patterns. In the hypercomplex domain, this redundancy strictly translates into a \textit{low-rank structure}.

Let $\Omega$ be the binary mask representing the spatial indices of the known (uncorrupted) pixels, and $\mathcal{P}_\Omega$ denote the orthogonal projection operator onto $\Omega$. The in-painting problem is traditionally formulated as a nuclear norm minimization problem, seeking to recover the lowest-rank matrix $Z$ that preserves the known visual data:
\begin{equation}
    \min_{Z} \|Z\|_* \quad \text{subject to} \quad \mathcal{P}_\Omega(Z) = \mathcal{P}_\Omega(Z_{orig})
\end{equation}

To efficiently solve this non-differentiable optimization problem, we implement a commutative hypercomplex adaptation of the SVT algorithm. This iterative technique leverages the SVD as a low-pass structural filter. By computing the SVD of the corrupted image $Z = U\Sigma V^H$\footnote{Let $A$ be an $\alpha\beta$-tessarine matrix. $A^H$ denotes its Hermitian transpose, defined as $A^{H_{3-2n}^1}$ where $n=1$ if $\alpha > 0$ and $n=2$ if $\alpha < 0$ (for further information see \cite{JimenezLopez2025}).} and applying a hard physical truncation at rank $k$, the algorithm discards the high-frequency components---which predominantly contain the null pixel mask (noise)---and forces the reconstruction of the missing gaps using the dominant eigenvectors. The iterative process is mathematically defined as follows:
\begin{enumerate}
    \item Compute the SVD of the current hypercomplex estimation: $Z^{(i)} = U \Sigma V^H$.
    \item Apply truncation to rank $k$ to extract the low-rank structural approximation: 
    \begin{equation}
        Z_{low} = \sum_{j=1}^{k} \sigma_j u_j v_j^H = U_k \Sigma_k V_k^H
    \end{equation}
    \item Enforce data consistency by updating only the missing regions, preserving original known pixels: 
    \begin{equation}
        Z^{(i+1)} = \mathcal{P}_{\Omega^C}(Z_{low}) + \mathcal{P}_\Omega(Z_{orig})
    \end{equation}
    where $\Omega^C$ is the complement of the mask.
\end{enumerate}

\subsubsection{Algorithmic Objective and Dimensional Bottlenecks}
The primary objective of this section is to benchmark the computational performance and scalability of the proposed \texttt{abtessarine\_Toolbox} (configured with $\alpha=-1$ and $\beta=1.5$) against the state-of-the-art quaternion algebra. 

While the \textit{Quaternion Toolbox for MATLAB} (QTFM), developed by Sangwine and Le Bihan \cite{Sanqwine_LeBihan}, serves as the foundational standard for hypercomplex signal processing, its native algebraic SVD decomposition is highly prohibitive for large-scale matrices. Therefore, to ensure a fair and competitive baseline, the comparative evaluation is not performed against the native QTFM SVD, but rather against a highly optimized complex-isomorphism truncated SVD implementation, denoted as \texttt{qsvds\_sp}, formulated by Gai et al. \cite{Gai2015}.

This optimized baseline maps the quaternion matrix $Q \in \mathbb{H}^{M \times N}$ into a symplectic complex matrix $M_c \in \mathbb{C}^{2M \times 2N}$ to leverage standard complex Arnoldi/Lanczos solvers (\texttt{svds}). However, this approach is inherently burdened by a \textit{dimensional inflation factor}. By expanding the problem space, the underlying solver must process a matrix four times larger in memory, carrying a severe $\mathcal{O}((2N)^3)$ computational complexity overhead. In contrast, the \texttt{abtessarine\_Toolbox} bypasses this bottleneck through the idempotent representation inherent to the commutative $\alpha\beta$-tessarine algebra. The hypercomplex SVD is mathematically decoupled into two independent $M \times N$ complex SVDs, strictly preserving the native spatial dimensions of the original image.

\subsubsection{Experiment I: Large-Scale Native Resolution Analysis}
\begin{figure*}[htbp]
    \centering
    \includegraphics[width=0.32\textwidth]{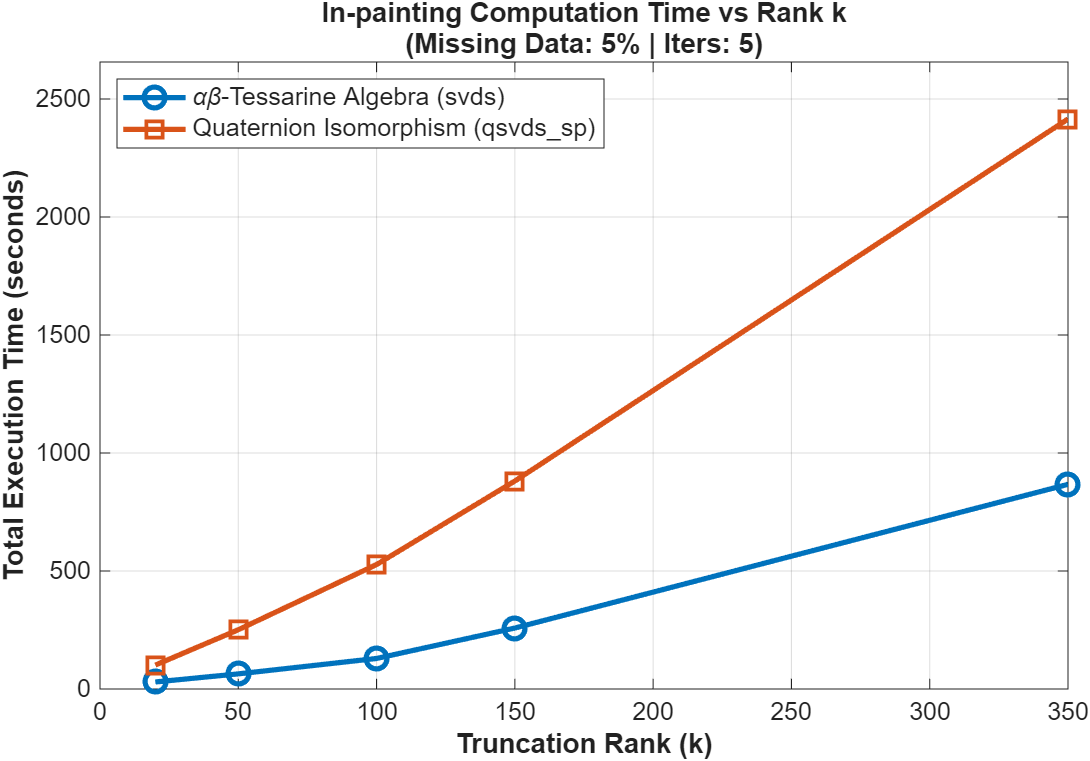}
    \hfill
    \includegraphics[width=0.32\textwidth]{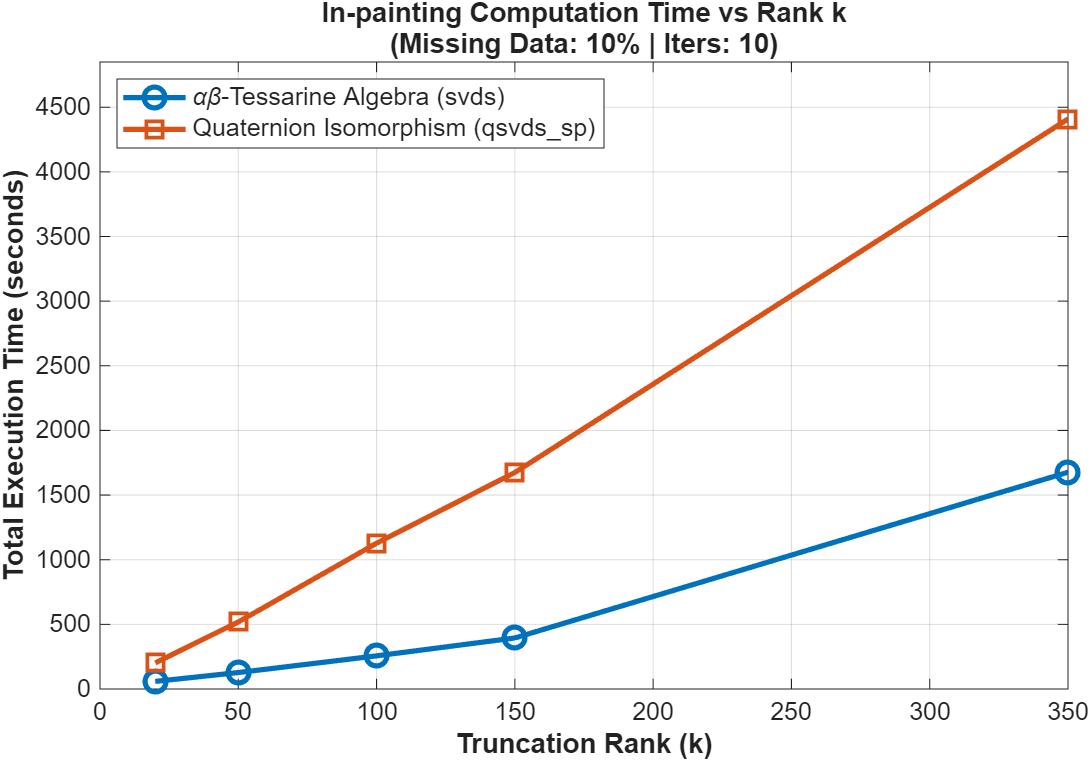}
    \hfill
    \includegraphics[width=0.32\textwidth]{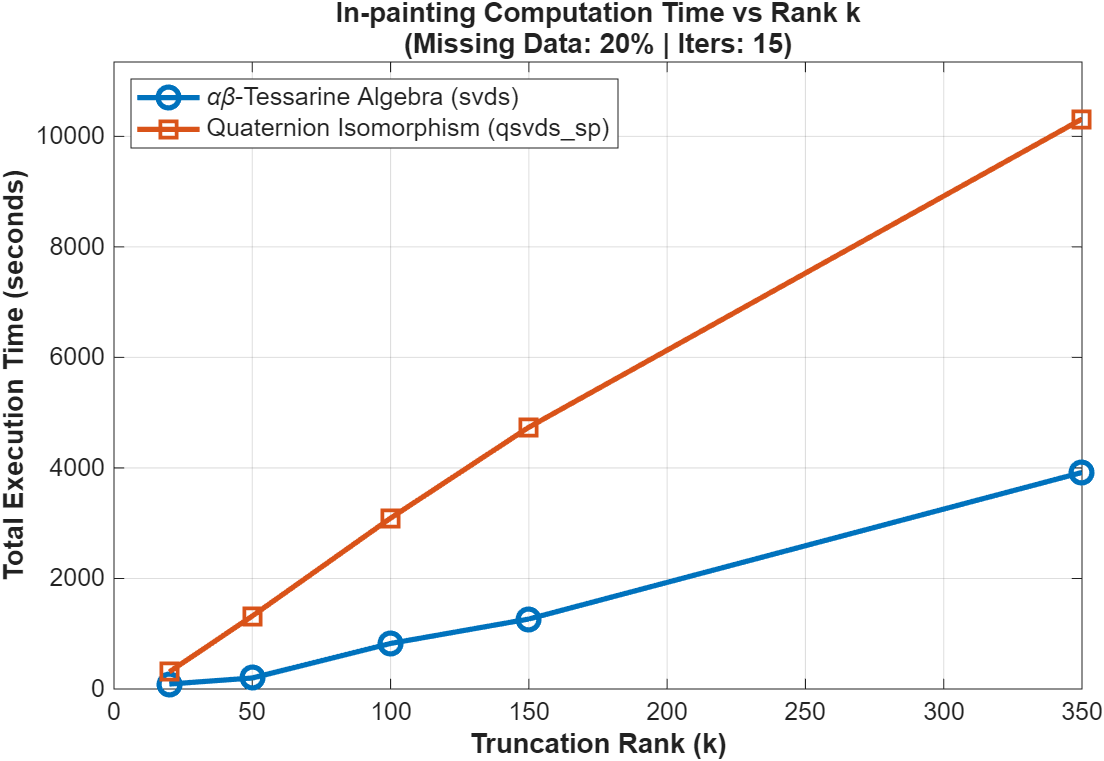}
    \caption{Computation time plots for the massive resolution dataset ($7952 \times 5304$) across 5\%, 10\%, and 20\% missing data. The quaternion baseline suffers from a severe exponential time explosion due to the $15904 \times 10608$ complex isomorphic mapping.}
    \label{fig:times_san_francisco}
\end{figure*}

To evaluate the impact of this dimensional inflation, the first benchmark was conducted using a massive high-resolution ``San Francisco'' image ($7952 \times 5304$ pixels). The reconstruction was tested under three missing data scenarios (5\%, 10\%, and 20\%), with adaptive iterations (5, 10, and 15, respectively) to prevent structural overfitting. The quantitative results are presented in Table~\ref{tab:inpainting_large} and visually supported by the execution time plots in Figure~\ref{fig:times_san_francisco}.

\begin{table*}[htbp]
\centering
\caption{Computational benchmark for color image in-painting at a massive native resolution ($7952 \times 5304$). Comparison between the proposed $\alpha\beta$-tessarine algebra and the optimized quaternion complex isomorphism.}
\label{tab:inpainting_large}
\renewcommand{\arraystretch}{1.2} 
\resizebox{\textwidth}{!}{%
\begin{tabular}{@{}lccccc@{}}
\toprule
\textbf{Rank ($k$)} & \multicolumn{2}{c}{\textbf{Execution Time (s)}} & \textbf{Speedup} & \multicolumn{2}{c}{\textbf{PSNR (dB)}} \\
\cmidrule(lr){2-3} \cmidrule(lr){5-6}
 & \textbf{$\boldsymbol{\alpha\beta}$-Tessarine} & \textbf{Quaternion} & & \textbf{$\boldsymbol{\alpha\beta}$-Tessarine} & \textbf{Quaternion} \\
\midrule
\multicolumn{6}{l}{\textit{Scenario 1: 5\% Missing Data, 5 Iterations}} \\
\midrule
20  & 29.79   & 101.34   & 3.40$\times$ & 34.31 & 34.27 \\
50  & 64.10   & 251.08   & 3.92$\times$ & 35.98 & 35.92 \\
100 & 129.42  & 527.39   & 4.07$\times$ & 37.54 & 37.44 \\
150 & 258.35  & 881.32   & 3.41$\times$ & 38.69 & 38.57 \\
350 & 867.23  & 2415.38  & 2.79$\times$ & 41.92 & 41.81 \\
\midrule
\multicolumn{6}{l}{\textit{Scenario 2: 10\% Missing Data, 10 Iterations}} \\
\midrule
20  & 58.79   & 202.31   & 3.44$\times$ & 31.29 & 31.25 \\
50  & 127.28  & 518.39   & 4.07$\times$ & 32.96 & 32.90 \\
100 & 256.31  & 1127.52  & 4.40$\times$ & 34.51 & 34.42 \\
150 & 394.72  & 1674.55  & 4.24$\times$ & 35.66 & 35.54 \\
350 & 1676.79 & 4408.86  & 2.63$\times$ & 38.85 & 38.74 \\
\midrule
\multicolumn{6}{l}{\textit{Scenario 3: 20\% Missing Data, 15 Iterations}} \\
\midrule
20  & 89.67   & 308.96   & 3.45$\times$ & 28.27 & 28.23 \\
50  & 199.22  & 1320.31  & 6.63$\times$ & 29.92 & 29.87 \\
100 & 823.88  & 3091.69  & 3.75$\times$ & 31.46 & 31.37 \\
150 & 1266.21 & 4736.59  & 3.74$\times$ & 32.58 & 32.48 \\
350 & 3918.65 & 10313.94 & 2.63$\times$ & 35.70 & 35.61 \\
\bottomrule
\end{tabular}%
}
\end{table*}

Regarding restoration fidelity, the PSNR values achieved by the \texttt{abtessarine\_Toolbox} are consistently superior to those computed by the quaternion baseline across all ranks. Although the differences are subtle (e.g., in Scenario 1, $41.92$ dB vs. $41.81$ dB at $k=350$), they confirm that the commutative algebra avoids the slight numerical degradation introduced by complex symplectic expansions. 

However, the most critical finding resides in the computational disparity. Focusing on Scenario 3, at high truncation ranks ($k=350$) and 20\% missing data, the quaternion approach requires over $10,313$ seconds to complete the iterative reconstruction, collapsing under the weight of its $15904 \times 10608$ mapped matrix. Conversely, the proposed toolbox completes the identical task in $3,918$ seconds, delivering sustained speedups of up to $6.63\times$.

\subsubsection{Experiment II: Algorithmic Stabilization on High-Frequency Textures}
\begin{figure*}[htbp]
    \centering
    \includegraphics[width=0.32\textwidth]{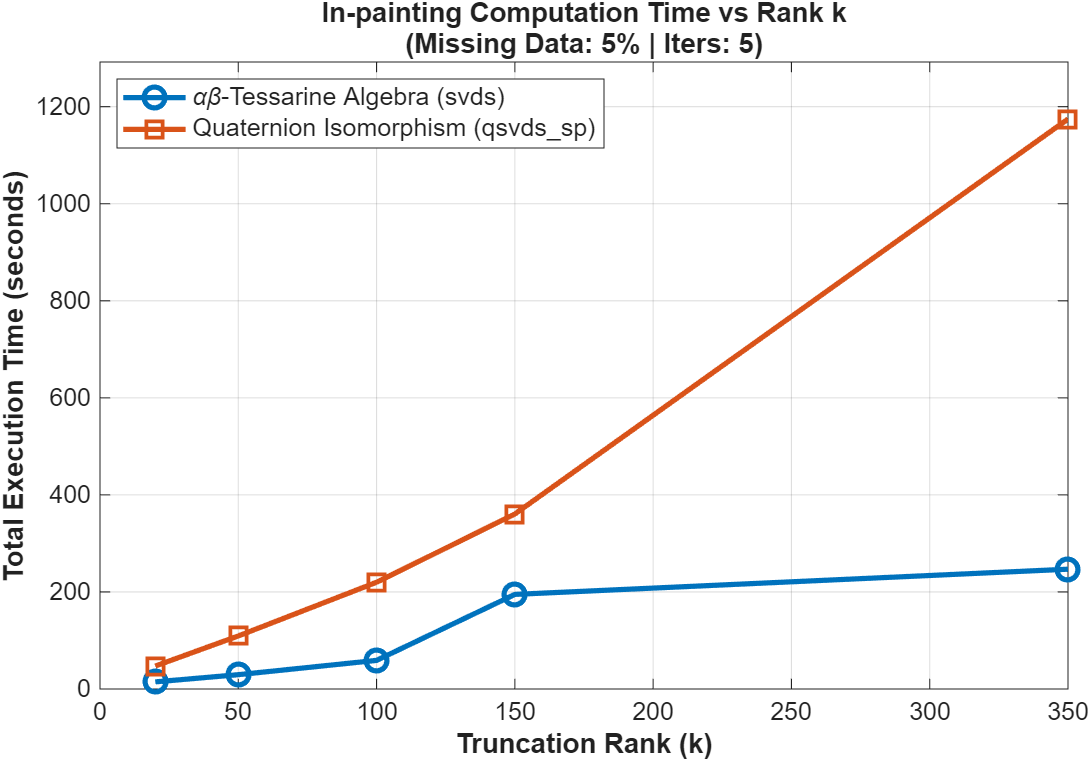}
    \hfill
    \includegraphics[width=0.32\textwidth]{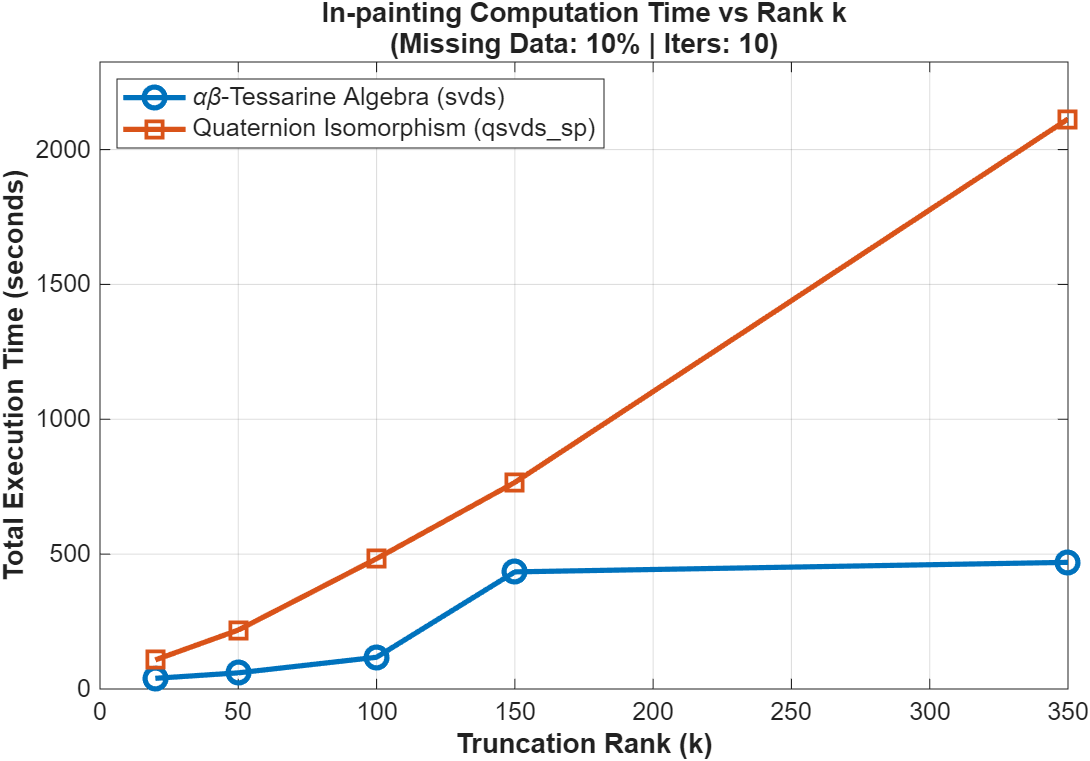}
    \hfill
    \includegraphics[width=0.32\textwidth]{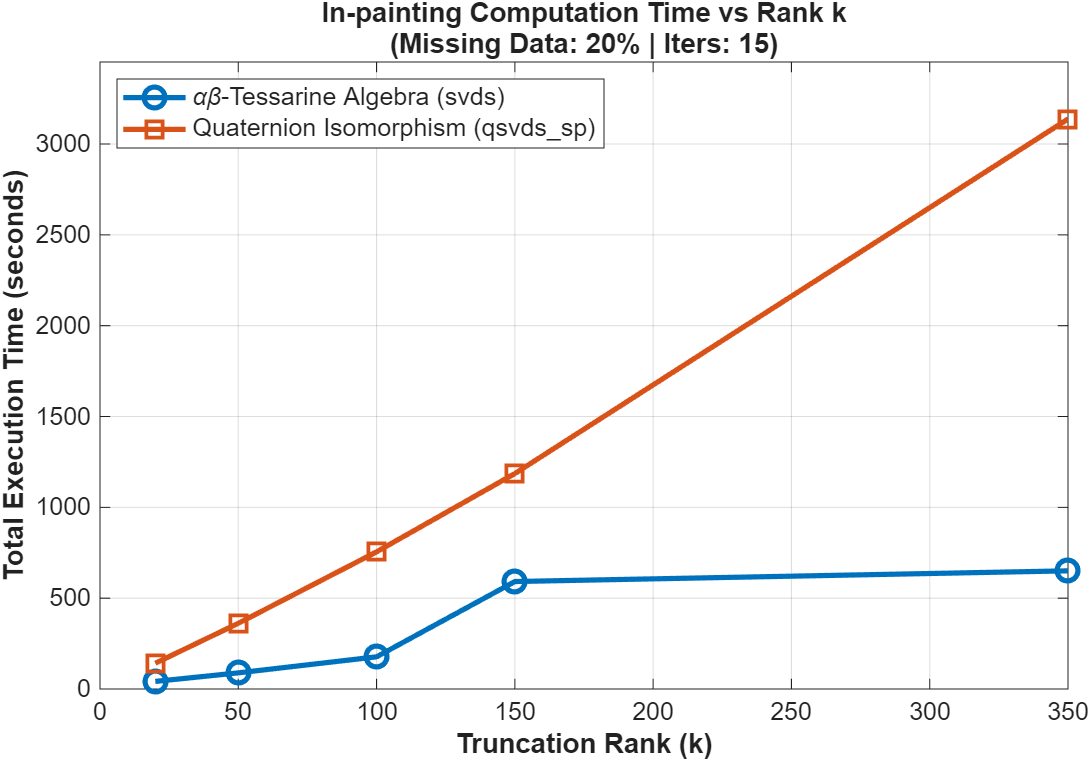}
    \caption{Computation time plots for the secondary dataset ($5464 \times 3640$). A pronounced algorithmic plateau is observed in the $\alpha\beta$-tessarine execution times at high truncation ranks ($k \ge 150$), demonstrating the ability of the proposed toolbox to safely transition to dense solvers and bound maximum computational times.}
    \label{fig:times_autumn_forest}
\end{figure*}

To rigorously analyze the behavior of the solver at extreme truncation ranks, a second experiment was conducted using an image with a spatial dimension of $3640 \times 5464$ but characterized by highly chaotic, high-frequency stochastic textures (``Autumn Forest''). The benchmark was evaluated under equivalent mathematical configurations (5\%, 10\%, and 20\% missing data), yielding the results presented in Table~\ref{tab:inpainting_forest} and Figure~\ref{fig:times_autumn_forest}.

\begin{table*}[htbp]
\centering
\caption{Computational benchmark for color image in-painting on a high-texture dataset ($5464 \times 3640$). The results highlight the algorithmic stabilization (plateau effect) of the $\alpha\beta$-tessarine SVD at high truncation ranks ($k \ge 150$).}
\label{tab:inpainting_forest}
\renewcommand{\arraystretch}{1.2}
\resizebox{\textwidth}{!}{%
\begin{tabular}{@{}lccccc@{}}
\toprule
\textbf{Rank ($k$)} & \multicolumn{2}{c}{\textbf{Execution Time (s)}} & \textbf{Speedup} & \multicolumn{2}{c}{\textbf{PSNR (dB)}} \\
\cmidrule(lr){2-3} \cmidrule(lr){5-6}
 & \textbf{$\boldsymbol{\alpha\beta}$-Tessarine} & \textbf{Quaternion} & & \textbf{$\boldsymbol{\alpha\beta}$-Tessarine} & \textbf{Quaternion} \\
\midrule
\multicolumn{6}{l}{\textit{Scenario 1: 5\% Missing Data, 5 Iterations}} \\
\midrule
20  & 14.67   & 47.33   & 3.23$\times$ & 33.16 & 33.11 \\
50  & 29.52   & 109.31  & 3.70$\times$ & 34.08 & 34.04 \\
100 & 58.99   & 219.98  & 3.73$\times$ & 34.87 & 34.82 \\
150 & 194.95  & 360.59  & 1.85$\times$ & 35.41 & 35.35 \\
350 & 246.87  & 1174.63 & 4.76$\times$ & 36.74 & 36.68 \\
\midrule
\multicolumn{6}{l}{\textit{Scenario 2: 10\% Missing Data, 10 Iterations}} \\
\midrule
20  & 39.58   & 107.85  & 2.72$\times$ & 30.14 & 30.10 \\
50  & 59.88   & 218.55  & 3.65$\times$ & 31.06 & 31.02 \\
100 & 117.73  & 483.19  & 4.10$\times$ & 31.84 & 31.80 \\
150 & 434.26  & 765.99  & 1.76$\times$ & 32.37 & 32.31 \\
350 & 469.58  & 2113.32 & 4.50$\times$ & 33.65 & 33.58 \\
\midrule
\multicolumn{6}{l}{\textit{Scenario 3: 20\% Missing Data, 15 Iterations}} \\
\midrule
20  & 41.78   & 141.90  & 3.40$\times$ & 27.13 & 27.08 \\
50  & 88.85   & 361.56  & 4.07$\times$ & 28.04 & 28.00 \\
100 & 177.09  & 755.11  & 4.26$\times$ & 28.80 & 28.75 \\
150 & 591.63  & 1186.04 & 2.00$\times$ & 29.29 & 29.24 \\
350 & 650.45  & 3137.71 & 4.82$\times$ & 30.45 & 30.40 \\
\bottomrule
\end{tabular}%
}
\end{table*}

As vividly illustrated in Figure~\ref{fig:times_autumn_forest}, a critical mathematical phenomenon emerges at higher ranks ($k \ge 150$). While the computation time for the quaternion baseline continues to grow relentlessly in an exponential fashion, the $\alpha\beta$-tessarine execution time experiences a drastic stabilization, effectively forming a computational plateau. For example, in Scenario 3, increasing the requested eigenvalues from $k=150$ to $k=350$ only increases the $\alpha\beta$-tessarine execution time by a marginal $58.82$ seconds (from $591.63$s to $650.45$s). By contrast, the quaternion baseline suffers a severe penalty of $1,951.67$ seconds for the exact same task.

This stabilization is a direct consequence of the internal heuristics employed by numerical SVD solvers. Truncated SVD algorithms (\texttt{svds}) rely on iterative Krylov subspace methods (e.g., ARPACK). As $k$ increases, the computational cost of orthogonalizing the Krylov basis vectors via Gram-Schmidt increases non-linearly. When $k$ crosses a specific heuristic threshold relative to the matrix dimensions, the solver determines that computing a full, dense SVD (via LAPACK) and manually truncating it is computationally cheaper than continuing the iterative approximation.

Because the \texttt{abtessarine\_Toolbox} operates strictly within the native $M \times N$ spatial dimensions, its matrices are compact enough to safely transition to the dense LAPACK solver. This effectively caps the maximum possible execution time, bounding it to a constant value regardless of how large $k$ becomes. In stark contrast, the dimensional inflation of the quaternion isomorphism ($7280 \times 10928$) creates a matrix so massive that executing a full dense SVD would cause a catastrophic Out-Of-Memory failure. Consequently, the quaternion algorithm is forced to remain trapped within the iterative ARPACK solver, suffering from severe orthogonalization penalties. This mathematical reality highlights that the commutative $\alpha\beta$-tessarine algebra not only provides substantial speedups but uniquely enables algorithmic safety and computational time-bounding in extreme-scale tensor processing.

\subsection{Scalability and Performance Evaluation: Dense SVD for Global Image Denoising}
\label{sec:denoising}

In this subsection, the raw algorithmic scalability and intrinsic mathematical performance of the \texttt{abtessarine\_Toolbox} are evaluated through a global image denoising framework. By applying a full, dense SVD to remove AWGN, the heuristics of iterative solvers are deliberately avoided in this stress test to expose the underlying computational complexity boundaries of the hypercomplex algebras. In the following benchmark, the execution times and restoration fidelity of the proposed commutative framework are rigorously compared against the optimized quaternion baseline, scaling spatial resolutions up to extreme dimensions to evaluate the critical hardware impact of the $\mathcal{O}(N^3)$ computational bottleneck.

\subsubsection{Mathematical Formulation of SVD-based Global Denoising}
To further evaluate the computational efficiency and mathematical robustness of the \texttt{abtessarine\_Toolbox}, a second benchmark is introduced focusing on global image denoising. Unlike the in-painting problem---which relies on iterative approximations to reconstruct localized missing data---this benchmark addresses the complete removal of AWGN across the entire spatial domain of an image.

Let $Z \in \mathbb{A}^{M \times N}$ be a hypercomplex image (where $\mathbb{A}$ is the chosen algebra) contaminated by AWGN with a variance $\sigma^2$. The pure mathematical approach to suppress this stochastic noise relies on the computation of a full, dense SVD:
\begin{equation}
    Z = U \Sigma V^H = \sum_{j=1}^{\min(M,N)} \sigma_j u_j v_j^H
\end{equation}

In this representation, the macroscopic structural components of the image (low frequencies) are strongly correlated and encoded within the largest singular values $\sigma_j$. Conversely, the uncorrelated Gaussian noise is uniformly distributed across the entire spectrum, dominating the smallest singular values. Therefore, a deterministic noise suppression is achieved by applying a hard physical truncation at a specific rank $k$:
\begin{equation}
    \hat{Z} = U_k \Sigma_k V_k^H \quad \text{where} \quad \Sigma_k = \text{diag}(\sigma_1, \dots, \sigma_k, 0, \dots, 0)
\end{equation}

This operation effectively annihilates the high-frequency stochastic noise, yielding a structurally denoised reconstruction $\hat{Z}$. From a computational perspective, the calculation of the full dense SVD represents one of the most demanding linear algebra operations, strictly bounded by an $\mathcal{O}(N^3)$ computational complexity for square matrices. This makes global denoising an ideal stress test to evaluate the intrinsic mathematical performance of hypercomplex representations, devoid of the heuristics present in iterative solvers.

\subsubsection{Experimental Setup and the ``Carina Nebula'' Benchmark}
To push the solvers to their computational limits, the benchmark utilizes the high-resolution ``Carina Nebula'' image provided by NASA. This specific dataset was selected due to its extreme high-frequency complexity; the dense distribution of stars, cosmic dust, and nebular gas makes it highly susceptible to over-smoothing (loss of fine detail) during the truncation process.

The evaluation compares the proposed \texttt{abtessarine\_Toolbox} (configured with $\alpha=-1$ and $\beta=1$) against the highly optimized quaternion complex-isomorphism SVD baseline (\texttt{qsvd\_sp}). To rigorously assess the $\mathcal{O}(N^3)$ scalability bottleneck, the image is progressively upsampled to square resolutions ranging from $256 \times 256$ up to an extreme $4096 \times 4096$ pixels. The experiment is conducted under three distinct AWGN damage scenarios (variance $\sigma^2 \in \{0.01, 0.05, 0.10\}$).

Crucially, to prevent structural degradation at higher scales, the truncation rank $k$ is logarithmically adapted to the spatial resolution (ranging from $k=20$ at $N=256$, to $k=200$ at $N=4096$). This adaptive strategy ensures that the essential high-frequency components (e.g., individual stars) are preserved as the pixel density increases, forcing the PSNR to scale positively with the resolution.

\subsubsection{Results and Computational Superiority}
\begin{figure*}[htbp]
    \centering
    \includegraphics[width=0.32\textwidth]{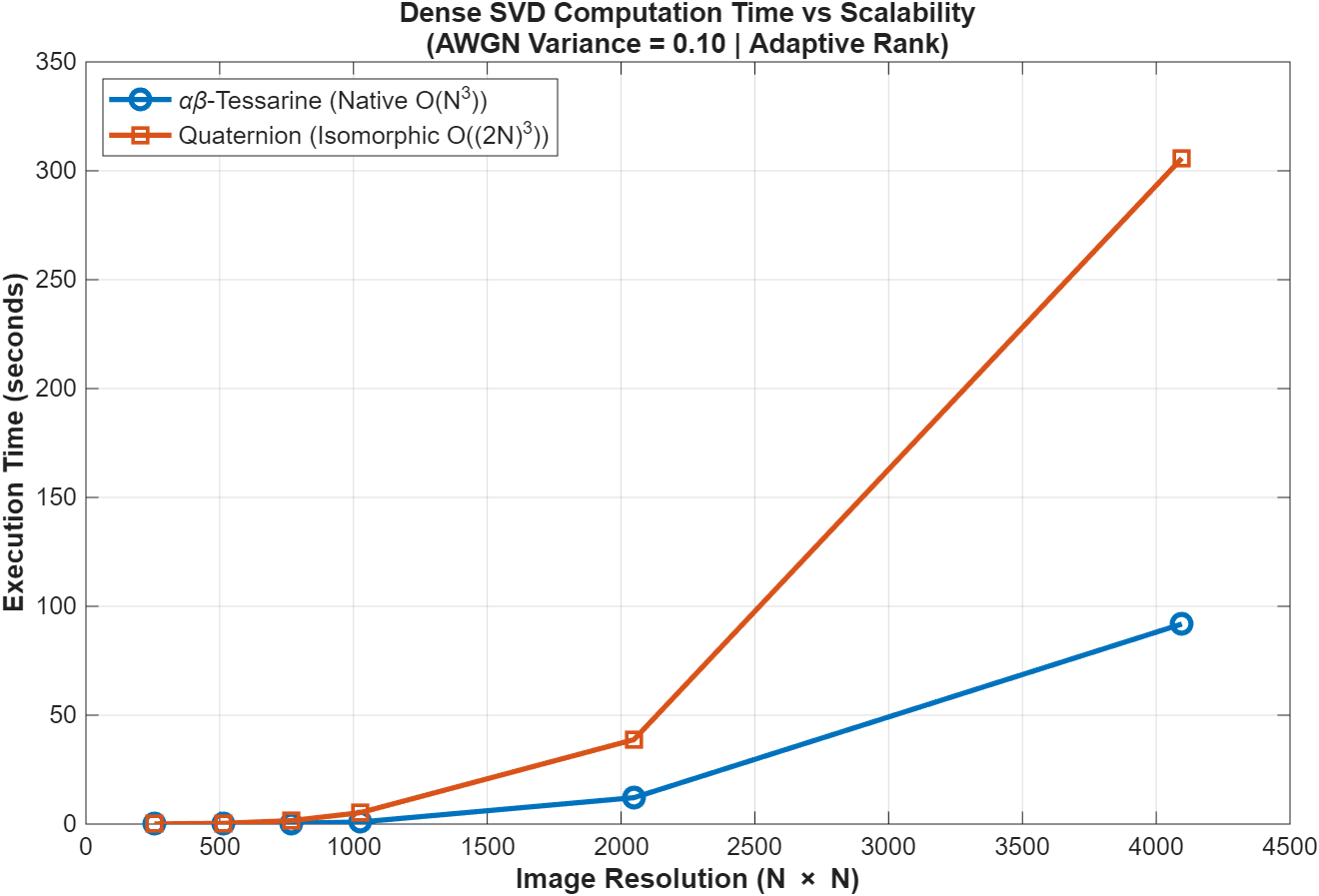}
    \hfill
    \includegraphics[width=0.32\textwidth]{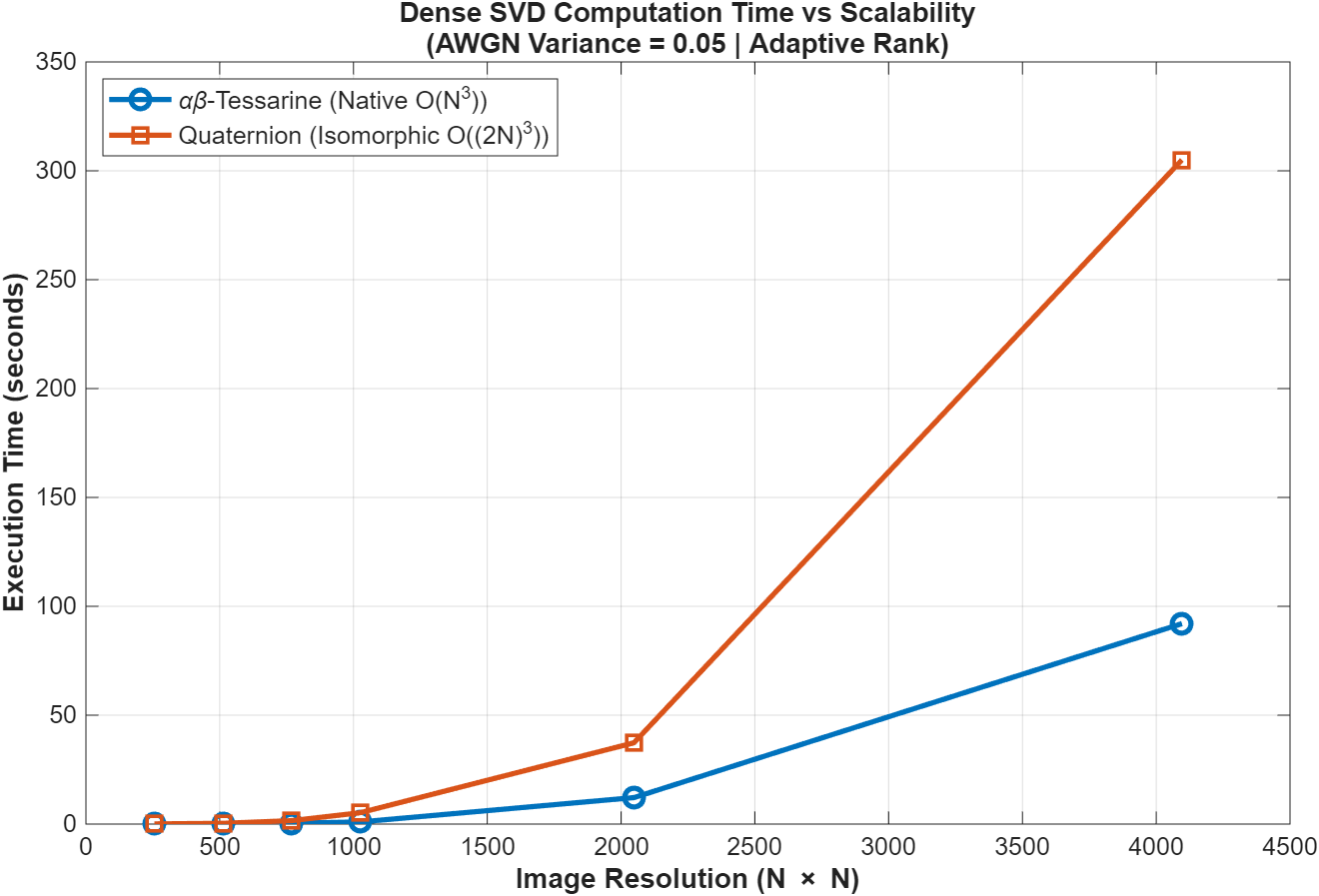}
    \hfill
    \includegraphics[width=0.32\textwidth]{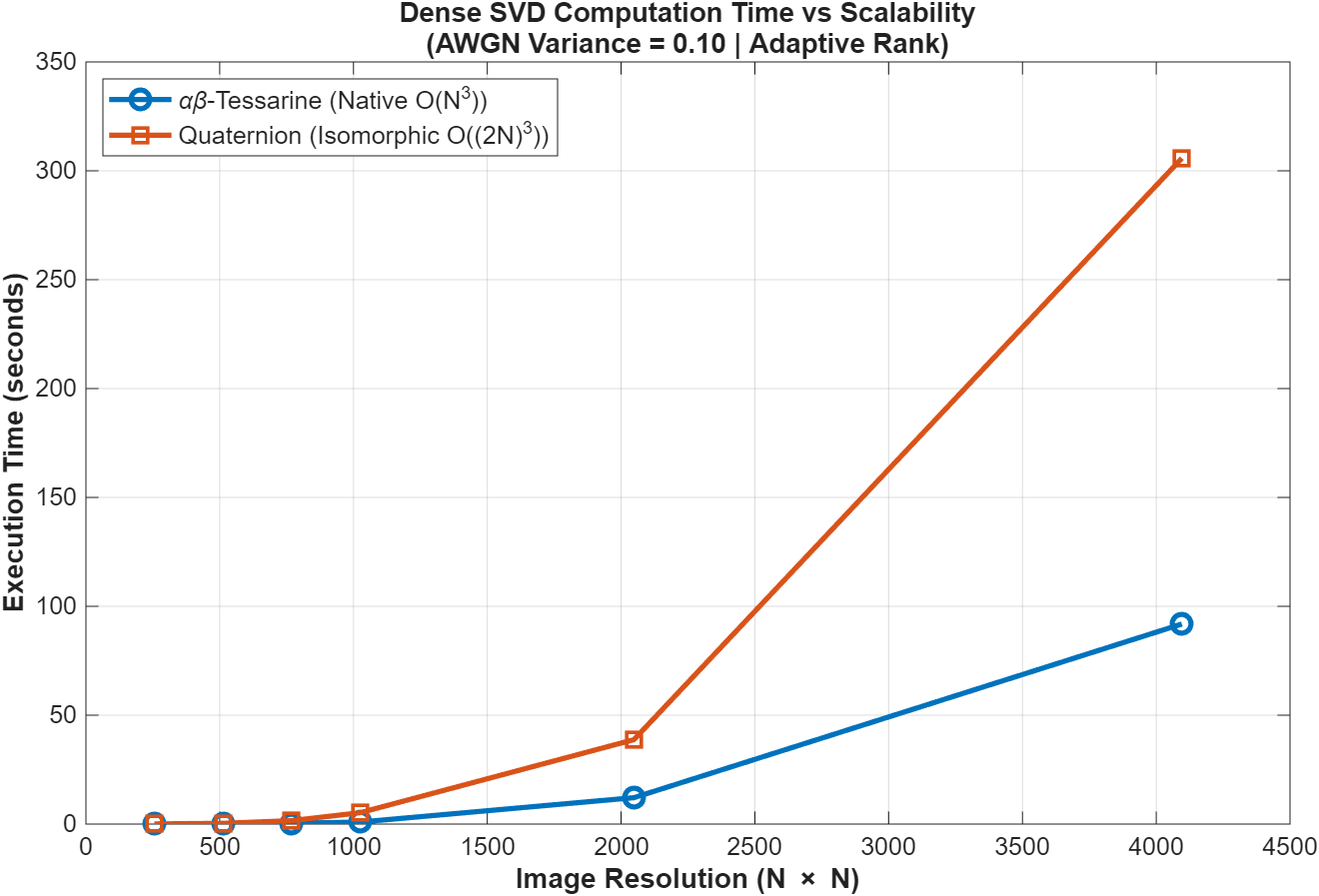}
    \caption{Computational scaling of the Dense SVD denoising algorithm across spatial resolutions. The native $\mathcal{O}(N^3)$ complexity of the \texttt{abtessarine\_Toolbox} exhibits a significantly shallower growth curve compared to the severe $\mathcal{O}((2N)^3)$ penalty induced by the quaternion isomorphism, effectively saving over 200 seconds of computation time at extreme scales ($4096 \times 4096$).}
    \label{fig:times_denoising}
\end{figure*}

The empirical results of the global denoising benchmark are detailed in Table~\ref{tab:denoising_benchmark} and visually supported by the computational time graphs in Figure~\ref{fig:times_denoising}.

\begin{table*}[htbp]
\centering
\caption{Computational benchmark for global image denoising (AWGN removal) using full Dense SVD. Comparison of execution times and restoration fidelity (PSNR) across progressively scaling native resolutions and varying noise variances.}
\label{tab:denoising_benchmark}
\renewcommand{\arraystretch}{1.2}
\resizebox{\textwidth}{!}{
\begin{tabular}{@{}lccrrcrr@{}}
\toprule
\textbf{Resolution} & \textbf{Rank} & \multicolumn{2}{c}{\textbf{Execution Time (s)}} & \textbf{Speedup} & \multicolumn{2}{c}{\textbf{PSNR (dB)}} \\
\cmidrule(lr){3-4} \cmidrule(lr){6-7}
\textbf{($N \times N$)} & \textbf{($k$)} & \textbf{Tessarine} & \textbf{Quaternion} & & \textbf{Tessarine} & \textbf{Quaternion} \\
\midrule
\multicolumn{7}{l}{\textit{Scenario 1: AWGN Variance = 0.01}} \\
\midrule
256 $\times$ 256   & 20  & 0.04  & 0.07   & 1.70$\times$ & 24.61 & 24.57 \\
512 $\times$ 512   & 30  & 0.16  & 0.35   & 2.25$\times$ & 24.82 & 24.79 \\
768 $\times$ 768   & 40  & 0.39  & 1.58   & 4.00$\times$ & 24.94 & 24.93 \\
1024 $\times$ 1024 & 50  & 0.94  & 5.25   & 5.58$\times$ & 25.08 & 25.06 \\
2048 $\times$ 2048 & 100 & 12.47 & 38.55  & 3.09$\times$ & 25.69 & 25.69 \\
4096 $\times$ 4096 & 200 & 92.58 & 308.76 & 3.34$\times$ & 26.67 & 26.67 \\
\midrule
\multicolumn{7}{l}{\textit{Scenario 2: AWGN Variance = 0.05}} \\
\midrule
256 $\times$ 256   & 20  & 0.06  & 0.08   & 1.50$\times$ & 19.75 & 19.70 \\
512 $\times$ 512   & 30  & 0.15  & 0.34   & 2.31$\times$ & 20.45 & 20.40 \\
768 $\times$ 768   & 40  & 0.40  & 1.57   & 3.93$\times$ & 20.77 & 20.69 \\
1024 $\times$ 1024 & 50  & 0.97  & 5.19   & 5.35$\times$ & 20.95 & 20.88 \\
2048 $\times$ 2048 & 100 & 12.15 & 37.33  & 3.07$\times$ & 21.14 & 21.06 \\
4096 $\times$ 4096 & 200 & 92.04 & 305.01 & 3.31$\times$ & 21.44 & 21.35 \\
\midrule
\multicolumn{7}{l}{\textit{Scenario 3: AWGN Variance = 0.10}} \\
\midrule
256 $\times$ 256   & 20  & 0.06  & 0.08   & 1.46$\times$ & 17.32 & 17.31 \\
512 $\times$ 512   & 30  & 0.15  & 0.34   & 2.27$\times$ & 18.03 & 17.99 \\
768 $\times$ 768   & 40  & 0.38  & 1.58   & 4.18$\times$ & 18.35 & 18.30 \\
1024 $\times$ 1024 & 50  & 0.93  & 5.21   & 5.62$\times$ & 18.52 & 18.47 \\
2048 $\times$ 2048 & 100 & 12.12 & 38.77  & 3.20$\times$ & 18.61 & 18.56 \\
4096 $\times$ 4096 & 200 & 91.87 & 305.79 & 3.33$\times$ & 18.78 & 18.73 \\
\bottomrule
\end{tabular}%
}
\end{table*}

The data robustly demonstrates the mathematical correctness of the adaptive truncation strategy: as the spatial resolution climbs, the PSNR progressively increases across all noise variances, confirming that the high-frequency cosmic structures are perfectly preserved while the un-correlated AWGN is mathematically filtered. Furthermore, the $\alpha\beta$-tessarine implementation consistently provides marginally superior PSNR values compared to the quaternion baseline, reinforcing the fidelity of the commutative algebraic decomposition.

However, the most critical finding resides in the computational execution times, where the limits of hardware architecture vividly expose the dimensional penalty of the quaternion formulation. Focusing on Scenario 1 (AWGN Variance = 0.01), at the $1024 \times 1024$ resolution, the \texttt{abtessarine\_Toolbox} achieves a peak speedup of $5.58\times$, processing the full SVD in under one second (0.94 seconds), while the baseline requires 5.25 seconds.

When the resolution is expanded to extreme scales ($4096 \times 4096$), the native $\alpha\beta$-tessarine solver maintains a highly efficient execution time of 92.58 seconds. In stark contrast, the quaternion algorithm collapses under the weight of its complex isomorphism. By expanding the data into an $8192 \times 8192$ complex matrix, the baseline SVD solver requires 308.76 seconds to complete a single decomposition. While the relative speedup naturally stabilizes at $3.34\times$ (due to the CPU transition from a compute-bound to a memory-bound state as the hardware cache saturates), the absolute time saved by the \texttt{abtessarine\_Toolbox} exceeds 216 seconds per image. This confirms that for modern, high-resolution tensor processing, preserving the native structural dimensions through commutative algebras is not merely an optimization, but a computational necessity.

\subsection{Illustrative Example: Digital Image Watermarking via Hypercomplex QR Decomposition}
\label{sec:watermarking}

In this final experimental subsection, the native QR decomposition capabilities of the \texttt{abtessarine\_Toolbox} are evaluated within the context of digital image watermarking. This application is specifically selected to assess the strict preservation of structural isometry, energy distribution, and numerical precision inherent to hypercomplex algebras. Through an energy-preserving spatial scalability benchmark, the execution times and restoration fidelity of the proposed commutative framework are rigorously compared against the state-of-the-art Structure-Preserving Quaternion QR (SP-QR) algorithm. This analysis is designed to demonstrate how the commutative idempotent representation completely circumvents the severe sequential bottlenecks and accumulated floating-point errors characteristic of iterative, non-commutative complex mapping algorithms.

\subsubsection{Mathematical Formulation and Isometry Requirements}
Digital image watermarking is a fundamental technique in multimedia security, utilized to embed copyright information or verification signals (the watermark) into a cover medium (the host image). To prevent the loss of inter-channel color correlations during the embedding process, hypercomplex algebras provide a natural framework to encapsulate the RGB channels into a single mathematical entity.

The QR decomposition serves as an exceptionally robust mathematical engine for this task. By decomposing a hypercomplex host image $Z_{host} \in \mathbb{A}^{M \times N}$ into an orthogonal matrix $Q$ and an upper triangular matrix $R$ ($Z_{host} = Q R$), the structural geometry of the image is separated from its energy distribution. The matrix $Q$ inherently preserves the directional spatial features, while the upper triangular matrix $R$ concentrates the signal's magnitude. Consequently, injecting the watermark into $R$ ensures that the hidden data is imperceptible yet resilient to external attacks.

The embedding process is defined as:
\begin{equation}
    R_{wm} = R + k_{adapt} Z_{wm}
\end{equation}
\begin{equation}
    Z_{marked} = Q R_{wm}
\end{equation}
where $Z_{wm}$ is the hypercomplex watermark and $k_{adapt}$ is the embedding gain. The embedded watermark can be blindly extracted for verification via $Z_{ext} = \frac{1}{k_{adapt}}(Q^H Z_{marked} - R)$.

\subsubsection{Algorithmic Baseline: Structure-Preserving QR}
To rigorously benchmark the numerical stability and computational performance of the proposed toolbox, the native $\alpha\beta$-tessarine QR routine was evaluated against the state-of-the-art Structure-Preserving Quaternion QR (SP-QR) algorithm proposed by Jia et al. \cite{Jia2018} .

Standard quaternion QR methods traditionally rely on complex adjoint mappings that often suffer from structural degradation, leading to a loss of strict orthogonality. The SP-QR algorithm brilliantly resolves this by employing symplectic Householder reflections on the $2N \times 2N$ complex isomorphic matrix, forcing the preservation of the quaternion structural constraints. 

However, this mathematical rigor comes at a severe computational cost. The SP-QR algorithm relies on an iterative, block-by-block processing scheme, executing thousands of sequential symplectic reflections. In contrast, the \texttt{abtessarine\_Toolbox} bypasses this iterative bottleneck entirely. By mapping the data to its idempotent representation, the computation is decoupled into independent commutative complex planes, allowing the direct execution of highly optimized LAPACK routines.

\subsubsection{Energy-Preserving Experimental Setup}
A spatial scalability benchmark was designed to stress-test the $\mathcal{O}(N^3)$ algorithmic complexity of the QR decomposition. The ``Autumn Forest'' high-resolution image was utilized as the host, with a ``Modern Architecture'' image acting as the watermark. The matrices were dynamically scaled to square resolutions ranging from $256 \times 256$ up to $2048 \times 2048$ pixels.

To conduct a scientifically rigorous evaluation of the PSNR across different scales, an \textit{Energy-Preserving Scaling} strategy was implemented. If a constant gain factor is applied indiscriminately as the matrix area expands quadratically, an disproportionate amount of energy is injected into the host image, leading to artificial degradation. Therefore, the embedding gain $k_{adapt}$ was mathematically scaled to preserve the total equivalent watermark energy relative to a baseline resolution of $N_{base} = 256$:
\begin{equation}
    k_{adapt} = k_{base} \left( \frac{N_{base}}{N} \right)
\end{equation}
The benchmark was executed under three specific scenarios, utilizing base gains of $k_{base} \in \{0.01, 0.05, 0.10\}$, representing light, medium, and aggressive embedding strengths, respectively.

\subsubsection{Results and Computational Scalability}
\begin{figure*}[htbp]
    \centering
    \includegraphics[width=0.32\textwidth]{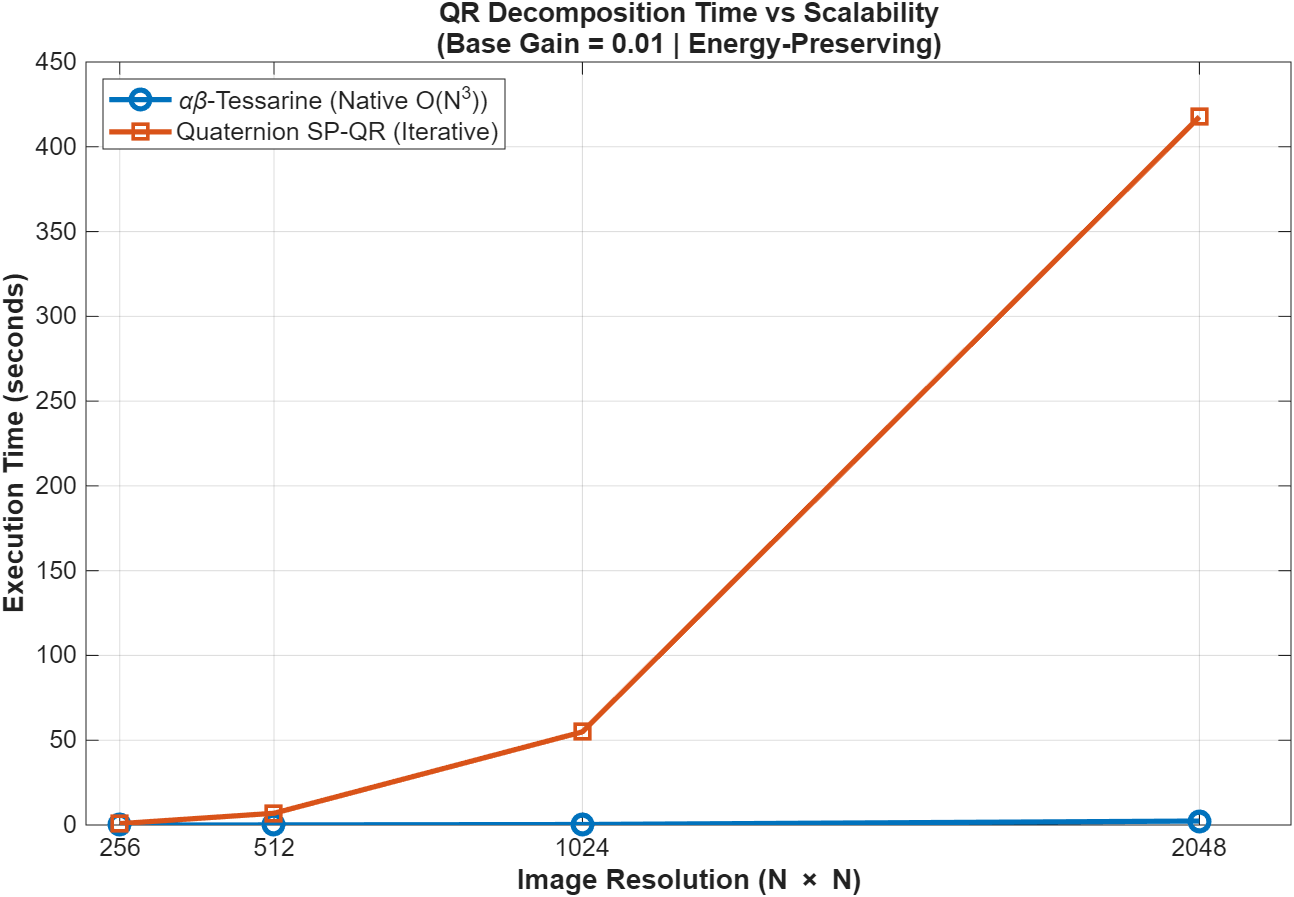}
    \hfill
    \includegraphics[width=0.32\textwidth]{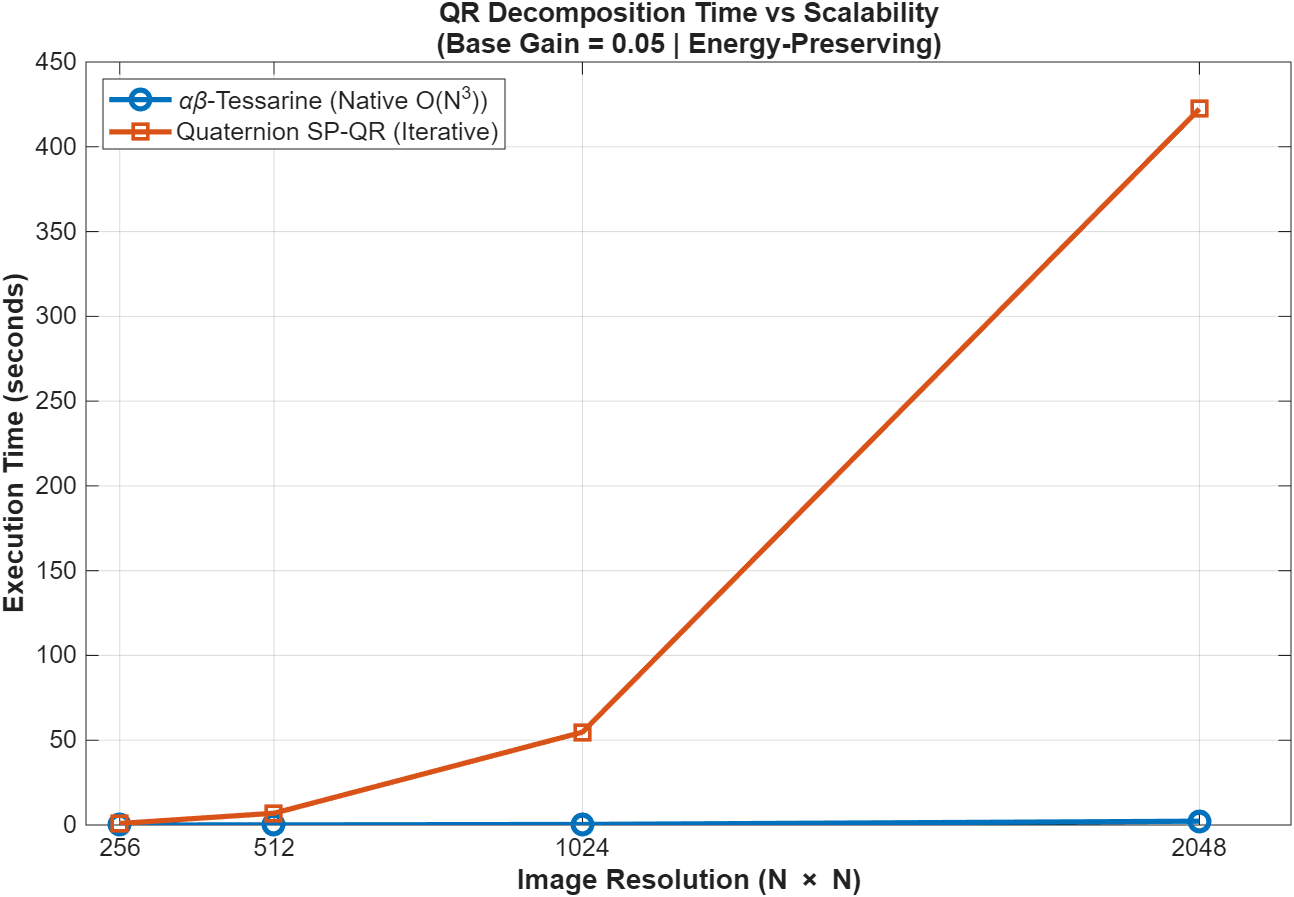}
    \hfill
    \includegraphics[width=0.32\textwidth]{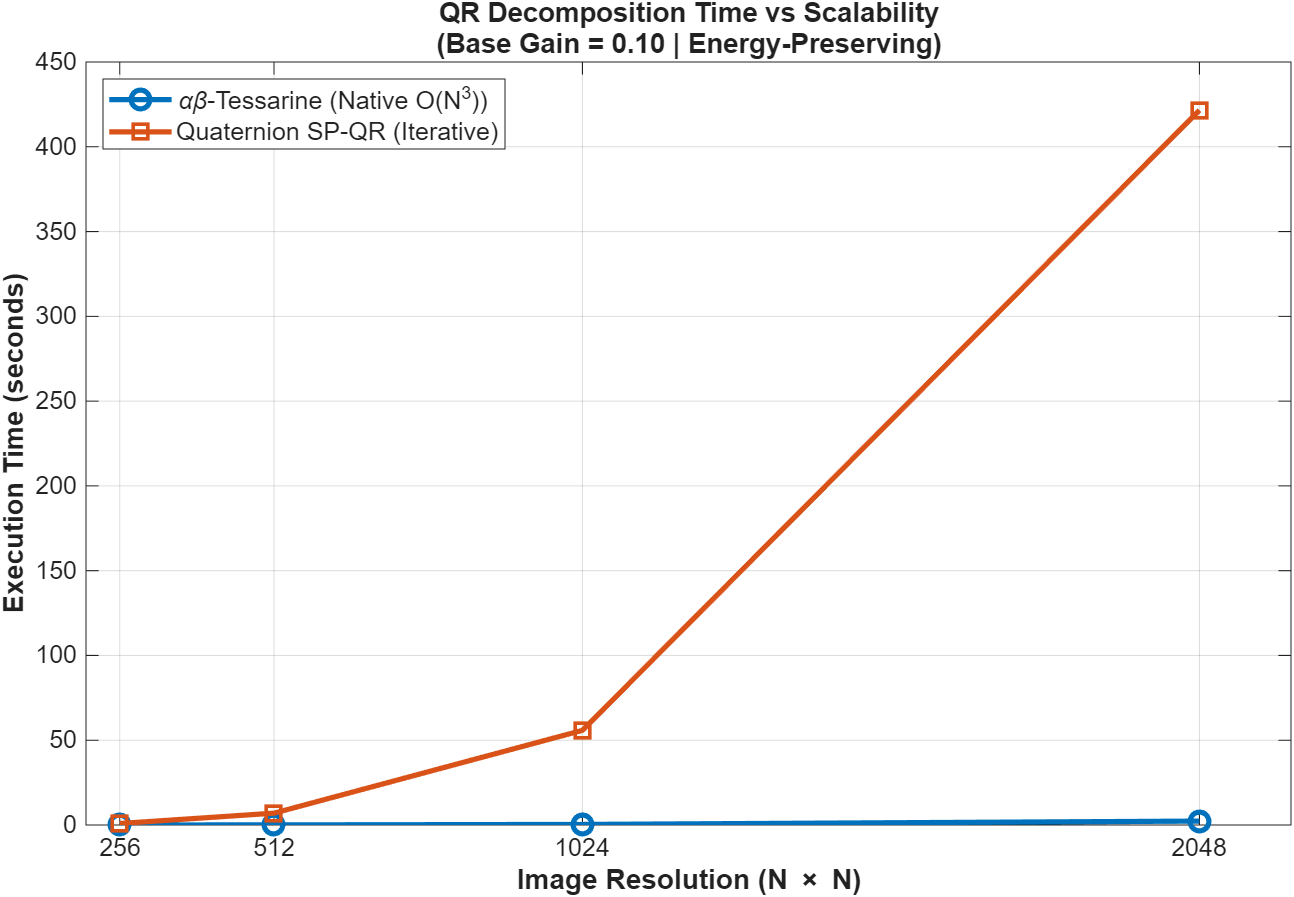}
    \caption{Computation time divergence in QR-based watermarking. The \texttt{abtessarine\_Toolbox} consistently processes $2048 \times 2048$ hypercomplex matrices in $\sim2.3$ seconds, bypassing the severe sequential bottleneck inherent to the structure-preserving iterative quaternion baseline.}
    \label{fig:times_watermarking}
\end{figure*}

The comprehensive results of the watermarking benchmark are detailed in Table~\ref{tab:watermarking_benchmark} and visually depicted in Figure~\ref{fig:times_watermarking}.

\begin{table*}[htbp]
\centering
\caption{Spatial scalability and energy-preserving benchmark for hypercomplex QR watermarking. Comparison of execution times and restoration fidelity between the \texttt{abtessarine\_Toolbox} and the iterative SP-QR baseline under varying embedding strengths.}
\label{tab:watermarking_benchmark}
\renewcommand{\arraystretch}{1.2}
\resizebox{\textwidth}{!}{%
\begin{tabular}{@{}llcrrcrr@{}}
\toprule
\textbf{Resolution} & \textbf{Adaptive Gain} & \multicolumn{2}{c}{\textbf{Execution Time (s)}} & \textbf{Speedup} & \multicolumn{2}{c}{\textbf{PSNR (dB)}} \\
\cmidrule(lr){3-4} \cmidrule(lr){6-7}
\textbf{($N \times N$)} & \textbf{($k_{adapt}$)} & \textbf{$\boldsymbol{\alpha\beta}$-Tessarine} & \textbf{Quaternion} & & \textbf{$\boldsymbol{\alpha\beta}$-Tessarine} & \textbf{Quaternion} \\
\midrule
\multicolumn{7}{l}{\textit{Scenario 1: Base Gain = 0.01 (Light Embedding)}} \\
\midrule
256 $\times$ 256   & 0.0100 & 0.0098 & 0.8034   & 82.31$\times$  & 48.45 & 47.06 \\
512 $\times$ 512   & 0.0050 & 0.0530 & 6.9232   & 130.71$\times$ & 54.37 & 51.80 \\
1024 $\times$ 1024 & 0.0025 & 0.2902 & 54.9284  & 189.26$\times$ & 60.06 & 55.97 \\
2048 $\times$ 2048 & 0.0013 & 2.3358 & 417.6321 & 178.79$\times$ & 64.99 & 59.49 \\
\midrule
\multicolumn{7}{l}{\textit{Scenario 2: Base Gain = 0.05 (Medium Embedding)}} \\
\midrule
256 $\times$ 256   & 0.0500 & 0.0366 & 0.8969   & 24.48$\times$  & 34.55 & 34.35 \\
512 $\times$ 512   & 0.0250 & 0.0615 & 6.9389   & 112.88$\times$ & 40.52 & 40.24 \\
1024 $\times$ 1024 & 0.0125 & 0.3338 & 54.7369  & 163.96$\times$ & 46.56 & 46.28 \\
2048 $\times$ 2048 & 0.0063 & 2.2468 & 422.2925 & 187.95$\times$ & 52.53 & 52.06 \\
\midrule
\multicolumn{7}{l}{\textit{Scenario 3: Base Gain = 0.10 (Aggressive Embedding)}} \\
\midrule
256 $\times$ 256   & 0.1000 & 0.0369 & 0.8816   & 23.90$\times$  & 28.61 & 28.41 \\
512 $\times$ 512   & 0.0500 & 0.0576 & 6.9807   & 121.09$\times$ & 34.54 & 34.34 \\
1024 $\times$ 1024 & 0.0250 & 0.3383 & 55.8294  & 165.01$\times$ & 40.58 & 40.50 \\
2048 $\times$ 2048 & 0.0125 & 2.3093 & 421.5448 & 182.54$\times$ & 46.58 & 46.56 \\
\bottomrule
\end{tabular}%
}
\end{table*}

The empirical data vividly illustrates the success of the energy-preserving mathematical model: across all base gains, the PSNR scales logarithmically upward as the spatial resolution increases. This validates that the structural integrity of the host image is perfectly maintained despite the quadratic growth of the matrix area. Furthermore, the $\alpha\beta$-tessarine algorithm consistently yields higher PSNR values than the SP-QR baseline across all scenarios. This discrepancy stems from the fundamental nature of the computation: while SP-QR strictly guarantees mathematical symplectic structure, its iterative execution of complex Householder reflections inevitably accumulates microscopic floating-point rounding errors. By leveraging the idempotent representation, the commutative $\alpha\beta$-tessarine algebra avoids iterative constraints, delivering exact numerical precision.

Crucially, the experiment exposes the profound computational bottleneck of non-commutative algorithmic constraints. The theoretical $\mathcal{O}(N^3)$ complexity of the QR decomposition dictates that a doubling in spatial resolution should result in roughly an eightfold increase in computation time. The quaternion SP-QR algorithm exhibits this exponential penalty dramatically, requiring  over 421 seconds in Scenario 3 to decompose a $2048 \times 2048$ image due to the $2N \times 2N$ complex expansion and the sequential nature of its Householder updates. 

Conversely, the \texttt{abtessarine\_Toolbox} absorbs this complexity with phenomenal efficiency, executing the identical matrix scale in approximately 2.3 seconds. This translates to an overwhelming, hardware-saturating speedup of up to $189\times$. The ability to process high-resolution hypercomplex decompositions in mere seconds rather than minutes demonstrates that the commutative, idempotent algebraic properties encapsulated within this toolbox are paramount for modern, real-time signal processing applications.

\section{Impact}

The utilization of hypercomplex algebras has profoundly impacted the fields of signal processing, computer vision, and multidimensional tensor algebra, among others. A prominent example of this scientific impact is the Quaternion Toolbox for MATLAB (QTFM) developed by Le Bihan et al. \cite{Sanqwine_LeBihan}, which has been extensively adopted by the scientific community. The availability of such software catalyzed a vast proliferation of research publications, proving that providing researchers with accessible hypercomplex tools directly accelerates scientific production.

More recently, the introduction of the toolbox for MATLAB about elliptic quaterniones \cite{Kosal2} has further expanded this research frontier. Its publication is currently driving a new wave of highly cited papers, demonstrating a strong, ongoing demand for computational frameworks based on these specific algebraic structures.

The development of the proposed \textit{$\alpha\beta$-Tessarine Toolbox} emerges as a direct response to the insistent demand from the scientific community for robust computational tools capable of solving real-world experimental problems using generalized hypercomplex models. By encompassing the generalized $\alpha\beta$-tessarine algebra, this toolbox implements a strictly more general mathematical framework than the algebra approached in \cite{Kosal2}. This inherently generalizes previous systems, enabling researchers to model and solve experimental problems that were previously constrained by the limitations of less generalized hypercomplex algebras.

Regarding the daily practice of its users, this software transitions the rigorous and often error-prone theoretical formulation of generalized hypercomplex numbers into an intuitive, high-performance computational environment. Unlike existing packages, this toolbox provides a substantially larger and more comprehensive suite of built-in scripts. This allows users to immediately tackle a wider array of applied experimental problems without writing algorithms from scratch, as demonstrated by the provided examples covering massive-scale color image in-painting via SVT, global image denoising via dense SVD, and energy-preserving digital image watermarking via QR decomposition.

Furthermore, by integrating high-performance computing capabilities, the toolbox opens new research questions in fields requiring massive scalability, such as large-scale machine learning and multidimensional data science \cite{Valle2021, Vieira2022, ZhangGao}, where computational bottlenecks previously rendered these generalized algebras impractical.

Finally, this toolbox is newly released as open-source software under the GPL-3.0 license via GitHub. While it has not yet generated spin-off companies or extensive usage metrics, its high-performance architecture is specifically designed to facilitate future integration into both large-scale academic research and commercial R\&D environments dealing with complex multidimensional signals.

\section{Conclusions}

In this work, the \textit{$\alpha\beta$-Tessarine Toolbox} is presented as a high-performance MATLAB framework for generalized hypercomplex tensor algebra. Through the implementation of the $\alpha\beta$-tessarine structure, a strictly more general mathematical environment is provided compared to existing state-of-the-art toolboxes, such as those dedicated to quaternions or specific tessarine subsets.

The gap between abstract hypercomplex theory and experimental application is successfully bridged by this toolbox. Due to its comprehensive suite of built-in scripts and high-performance architecture, multidimensional problems---ranging from massive-scale color image in-painting and global image denoising to energy-preserving digital watermarking---can be tackled with unprecedented computational efficiency.

It is concluded that the \textit{$\alpha\beta$-Tessarine Toolbox} represents a significant advancement in the availability of computational tools for the scientific community. Its open-source nature and scalability lay the groundwork for future research in large-scale machine learning and multidimensional signal processing, where new insights and superior performance over traditional methods can be offered by the generalized nature of $\alpha\beta$-tessarine algebra.
\section*{Acknowledgements}
\label{}
This work was supported by the Research and Knowledge Transfer Plan
2025 of the University of Ja\'en (grant 2025/00345/001).

 \appendix

\section{Cholesky Factorization in the $(\alpha\beta)$-Tessarine Algebra}\label{Cholesky factorization}

Given a hypercomplex Hermitian positive-definite matrix $X \in \mathbb{T}_{\alpha\beta}^{n \times n}$, its factorization is sought in the form $X = LL^{H_{3-2n}^1}$, where $n=1$ if $\alpha > 0$ and $n=2$ if $\alpha < 0$. Here, $L$ is a hypercomplex lower triangular matrix and $L^{H_{3-2n}^1}$ represents its Hermitian transpose.

Using the idempotent projector basis of the algebra, $w_1 = \frac{\sqrt{\beta}+j}{2\sqrt{\beta}}$ and $w_2 = \frac{\sqrt{\beta}-j}{2\sqrt{\beta}}$, the matrices are decomposed into their continuous branches:
\[
    X = X_s w_1 + X_d w_2
\]
\[
    L = L_s w_1 + L_d w_2
\]
\[
    L^{H_{3-2n}^1} = L_s^{H_{3-2n}^1} w_1 + L_d^{H_{3-2n}^1} w_2
\]

Substituting into $X = LL^{H_{3-2n}^1}$ and applying the orthogonality of the basis ($w_1 w_2 = 0$, $w_1^2 = w_1$, $w_2^2 = w_2$), the problem is decoupled into two independent Cholesky factorizations in the $(\alpha)$-complex domain:
\[
    X_s = L_s L_s^{H_{3-2n}^1} \quad \text{and} \quad X_d = L_d L_d^{H_{3-2n}^1}
\]

The branches $X_s$ and $X_d$ belong to the $(\alpha)$-complex subspace. The analytical and computational resolution requires projecting these matrices into standard domains depending on the sign of $\alpha$.
\begin{itemize}
\item For $\alpha < 0$ (where $n=2$), the matrices $X_s = A_s + B_s i$ and $X_d = A_d + B_d i$ operate with $i^2 = \alpha$. These matrices are projected into the standard complex field $\mathbb{C}$ using the classical imaginary unit $\epsilon$ ($\epsilon^2 = -1$).

The projected matrices $\dot{X}_s, \dot{X}_d \in \mathbb{C}^{n \times n}$ are defined by scaling the imaginary parts by $\sqrt{|\alpha|}$:
\[
    \dot{X}_s = A_s + \dot{B}_s \epsilon, \quad \text{where} \quad \dot{B}_s = \sqrt{|\alpha|}B_s
\]
\[
    \dot{X}_d = A_d + \dot{B}_d \epsilon, \quad \text{where} \quad \dot{B}_d = \sqrt{|\alpha|}B_d
\]

Since $X$ is 2-Hermitian and positive-definite, the projections $\dot{X}_s$ and $\dot{X}_d$ are classical complex Hermitian positive-definite matrices. The standard Cholesky factorization is applied to obtain the lower triangular factors:
\[
    \dot{L}_s = \text{chol}(\dot{X}_s, \text{'lower'})
\]
\[
    \dot{L}_d = \text{chol}(\dot{X}_d, \text{'lower'})
\]

To return to the $(\alpha)$-complex subspace, the real parts are retained and the imaginary parts are divided by the topological factor $\sqrt{|\alpha|}$:
\[
    L_s = \Re(\dot{L}_s) + \frac{\Im(\dot{L}_s)}{\sqrt{|\alpha|}} i
\]
\[
    L_d = \Re(\dot{L}_d) + \frac{\Im(\dot{L}_d)}{\sqrt{|\alpha|}} i
\]

\item For $\alpha > 0$ (where $n=1$), the Minkowski topology requires projecting the branches into purely real, symmetric, positive-definite matrices:
\[
    E_{s1} = A_s + \sqrt{\alpha} B_s, \quad E_{s2} = A_s - \sqrt{\alpha} B_s
\]
\[
    E_{d1} = A_d + \sqrt{\alpha} B_d, \quad E_{d2} = A_d - \sqrt{\alpha} B_d
\]

Four purely real Cholesky factorizations are executed:
\[
    L_{s1} = \text{chol}(E_{s1}, \text{'lower'}), \quad L_{s2} = \text{chol}(E_{s2}, \text{'lower'})
\]
\[
    L_{d1} = \text{chol}(E_{d1}, \text{'lower'}), \quad L_{d2} = \text{chol}(E_{d2}, \text{'lower'})
\]

The $(\alpha)$-complex triangular components are reconstructed by recombining the real factors:
\[
    L_s = \frac{L_{s1} + L_{s2}}{2} + \frac{L_{s1} - L_{s2}}{2\sqrt{\alpha}} i
\]
\[
    L_d = \frac{L_{d1} + L_{d2}}{2} + \frac{L_{d1} - L_{d2}}{2\sqrt{\alpha}} i
\]
\end{itemize}

Finally, the hypercomplex lower triangular matrix is uniformly reconstructed in the algebra by substituting back into the idempotent basis:
\[
    L = \frac{L_s + L_d}{2} + \frac{L_s - L_d}{2\sqrt{\beta}} j
\]

\section*{Current executable software version}
\label{sec:executable}

Not applicable. The \textit{$\alpha\beta$-Tessarine Toolbox} is distributed entirely as uncompiled MATLAB source code scripts (\texttt{.m} files). It does not require a standalone executable, as it is executed directly within the native MATLAB environment. The complete source code is accessible via the repository provided in Table \ref{code_metadata}.


\begin{thebibliography}{00}

\bibitem{Hahn2016}
S.L.  Hahn,  K.M.  Snopek, {\em Complex and Hypercomplex Analytic Signals: Theory and Applications},  Artech House: Norwood, MA, USA (2016).

\bibitem{Valle2021}
 M.E. Valle,   R.A. Lobo,  Hypercomplex-valued recurrent correlation neural networks,  {\em Neurocomputing} 432 (2021)  111--123.


\bibitem{Vieira2022}
G. Vieira,   M.E. Valle, A general framework for hypercomplex-valued extreme learning machines,  {\em J.  Comput.  Math.  Data Sci.} 3 (2022) 100032.   

\bibitem{ZengSong}
R. Zeng, Q. Song, S. Sun, Quaternion version of the It\^o's formula, {\em Math. Comput. Simul.} 222 (2024) 242--251.

\bibitem{Borio}
D. Borio, M. Susi, GNSS meta-signal tracking using a bicomplex Kalman filter,  {\em Navigation-US}  71~(4) (2024) 674.


\bibitem{ZhangGao}
H. Zhang, P. Gao, R. Ye, I. Stamova, J. Cao, Fixed/Predefined time synchronization of fractional quaternion delayed neural networks with disturbances, {\em Math. Comput. Simul.} 232 (2025) 276--294.


\bibitem{Pei}
S.C. Pei, J.H. Chan, J.J. Ding,  M.Y. Chen,  Eingenvalues and singular value descompostions of reduced biquaternion matrices,  {\em IEEE Trans. Circuits Syst.}  55~(9) (2008)	2673--2685.

\bibitem{Ortolani2017}
 F.  Ortolani,  M. Scarpintini,  D. Comminiello,  A. Uncini,  On 4-dimensional hypercomplex algebras in adaptive signal processing,  in: A. Esposito, M. Faundez-Zanuy, F. Morabito,  E. Pasero (Eds.), {\em Neural Adv. Process. Nonlinear Dyn. Signals},  Springer (2017) 131--140.  

\bibitem{Cariow2024}
A.  Cariow,  G. Cariowa,   Reduced-complexity algorithms for tessarine neural networks,  {\em IEEE Trans.   Neural  Netw.   Learn.  Syst.}  (2024) 1--6.  

\bibitem{Navarro_beta_quaternion}
J.  Navarro-Moreno, R.M.  Fern\'andez-Alcal\'a, J.D.  Jim\'enez-L\'opez, J.C. Ruiz-Molina,  Proper adaptive filtering in four dimensional Cayley-Dickson algebras,  {\em J. Frank. Inst.}   360~(12) (2023)  7739--7769.  

\bibitem{Navarro_Segre}
J. Navarro-Moreno,    R.M. Fernández-Alcalá,   J.C. Ruiz-Molina,  Proper ARMA modelling and forecasting in the generalized Segre’s quaternion domain,  {\em Mathematics}  10~(7) (2022) 1083.

\bibitem{Kosal1}
H.H. K{\"{o}}sal, E. Ki{\c{s}}i, M. Akyi{\u{g}}it, B. {\c{C}}elik, 
Elliptic quaternion matrices: Theory and algorithms, 
{\em Axioms} 13(10) (2024) 656.

\bibitem{Kosal3}
G. Atali, H.H. K{\"{o}}sal, M. Pekyaman, 
A new image restoration model associated with special elliptic quaternionic least-squares solutions based on LabVIEW, 
\textit{J. Comput. Appl. Math.} 425 (2023) 115071. 

\bibitem{Sanqwine_LeBihan}
S.J. Sangwine, N. Le Bihan, 
Quaternion Toolbox for MATLAB (QTFM), 
\textit{Software package} (2005). 
https://qtfm.sourceforge.io/

\bibitem{Kosal2}
H.H.  K\"osal,  E.  Ki\c{s}i,  M. Akyi\u{g}it,  B. \c{C}elik,  Elliptic quaternion matrices: A MATLAB toolbox and applications for image processing, {\em Axioms} 13~(11) (2024) 771.

\bibitem{JimenezLopez2025}
J.D. Jim{\'{e}}nez-L{\'{o}}pez, J. Navarro-Moreno, R.M. Fern{\'{a}}ndez-Alcal{\'{a}}, J.C. Ruiz-Molina, 
Advancing Computational Tools for Analyzing Commutative Hypercomplex Algebras, 
\textit{arXiv preprint arXiv:2508.02709} (2025). 

\bibitem{Bihan2017}
N. Le Bihan, 
The geometry of proper quaternion random variables, 
\textit{Signal Process.} 138 (2017) 106--116.

\bibitem{Nitta2019}
T. Nitta, M. Kobayashi, D.P. Mandic, 
Hypercomplex widely linear estimation through the lens of underpinning geometry, 
\textit{IEEE Trans. Signal Process.} 67 (2019) 3985--3994.

\bibitem{Grassucci_proper}
E. Grassucci, D. Comminiello, A. Uncini, 
An information-theoretic perspective on proper quaternion variational autoencoders, 
\textit{Entropy} 23(7) (2021) 856.

\bibitem{Gai2015}
S. Gai, G. Yang, M. Wan, L. Wang, Denoising color images by reduced quaternion matrix singular value decomposition, {\em Multidimens. Syst. Signal Process.} 26 (2015) 307–320.



\bibitem{Jia2018}
Z.G. Jia, M.S. Wei, M.X. Zhao, Y. Chen, 
A new real structure-preserving quaternion QR algorithm, 
{\em J. Comput. Appl. Math.} 343 (2018) 26--48. 

\end{thebibliography}
\end{document}